\documentclass[12pt]{amsart}

\usepackage{graphicx}
\usepackage{amsxtra}
\usepackage{amssymb}
\usepackage{amscd}
\usepackage{longtable}
\usepackage{array}
\usepackage[pdftex]{hyperref}

\newcommand\id{\operatorname{Id}}

\newcommand\rank{\operatorname{rank}}

\theoremstyle{plain}

\newtheorem{thm}{Theorem}[section]

\newtheorem{prop}[thm]{Proposition}
\newtheorem{cor}[thm]{Corollary}

\newtheoremstyle{osa}{6pt}{6pt}{\rmfamily}{0pt}{\rmfamily}{.}{1 ex}{}

\theoremstyle{osa}

\newtheorem{ex}[thm]{Example}

\newtheorem{remark}{Remark}

\title[Hermitian structures on cotangent Lie groups]{Hermitian structures on
cotangent bundles\\ of four dimensional solvable Lie groups}

\author[L.C. de Andr\'es, M. L. Barberis, I. Dotti, M. Fern\'andez]{Luis C. de
Andr\'es, M. Laura Barberis, Isabel Dotti, Marisa Fern\'andez}

\date{
}

\begin{document}

\begin{abstract}
We study  hermitian structures, with respect to the standard
neutral metric on the cotangent bundle $T^*G$ of a 2n-dimensional
Lie group $G$, which are left invariant with respect to the Lie
group structure on $T^*G$ induced by the coadjoint action. These
are in one-to-one correspondence with  left invariant generalized
complex structures on $G$. Using this correspondence and results
of \cite{CG} and \cite{FGG}, it turns out that when $G$ is
nilpotent and four or six dimensional, the cotangent bundle $T^*G$
always has a hermitian structure. However, we prove that if $G$ is
a four dimensional  solvable
Lie group admitting neither complex
nor symplectic structures, then $T^*G$ has no hermitian structure or,
equivalently, $G$ has no left invariant generalized complex structure.

\end{abstract}

\maketitle

\section{Introduction}
The cotangent bundle $T^*G$ of a Lie group $G$ with Lie algebra
$\mathfrak g$ has a canonical Lie group structure induced by the
coadjoint action of $G$ on $\mathfrak g^*$ and also a canonical
bi-invariant neutral metric. With respect to this data,  hermitian
structures on $T^*G$ such that left translations are holomorphic
isometries are given by endomorphisms $J$ of $\mathfrak g \oplus
\mathfrak g^*$ satisfying $J^2=-\id$ which are orthogonal with
respect to
\begin{equation} \label{st-bilform}
\langle (x,\alpha), (y, \beta)\rangle = \frac12 \left(\beta(x)+\alpha(y)\right),
\end{equation}
and satisfy $N_J\equiv 0$, where $N_J$ is
defined in \eqref{nijen},
with respect to
the Lie bracket:
\begin{equation} \label{bracket1}
[ (x,\alpha), (y, \beta)] = ([x,y],-\beta \circ \mbox{ad}(x) +
\alpha \circ \mbox{ad}(y)) \hspace{1cm} \mbox{for }x, y\in \mathfrak
g, \;\; \alpha,\; \beta \in  \mathfrak g^* .
\end{equation}
On the other hand, $\mathfrak g \oplus \mathfrak g^*$ is the fiber
at the identity $e$ of the bundle $TG \oplus T^*G$ over $G$ and
one may extend $J$ above to the whole $TG \oplus T^*G$ using the
standard lift of left multiplication in $G$.  The Courant bracket
(see \eqref{courant} below), when restricted to left invariant
vector fields and left invariant 1-forms is given by the equation
above thus establishing a correspondence, in the invariant case,
between invariant hermitian structures on $T^*G$ and left
invariant generalized complex structures on $G$
(Proposition~\ref{corresp}). It follows that any such structure
gives rise to a Poisson Lie group structure on $T^*G$ such  that
the dual Poisson Lie group $(T^*G)^*$ is a complex Lie group
(Corollary \ref{Poisson}).

The concept of \emph{generalized complex structure} was introduced
by Hitchin~\cite{Hitchin} and developed by Gualtieri~\cite{Gualt}.
Symplectic and complex geometry are extremal special cases of
generalized complex geometry. In~\cite{CG} Cavalcanti and
Gualtieri show that the $34$ classes of $6$--dimensional nilpotent
Lie groups (see~\cite{Mag,Sal} for the classification) have a left
invariant generalized complex structure; but, five of these
classes of nilpotent Lie groups admit neither symplectic nor
complex left invariant geometries (see~\cite{Sal}). It is proved
in \cite{FGG}  that every four dimensional nilpotent Lie group has
left invariant symplectic structures and hence generalized complex
structures. So, it seems interesting to understand the way this
property occurs on non-nilpotent solvable Lie groups.

In this paper we deal with left invariant generalized complex
structures on  solvable Lie groups of dimension $4$. To this end,
in Proposition \ref{corresp} of \S\ref{cotg}, we show that there
is a one-to-one correspondence between left invariant generalized
complex  structures on a Lie group $G$ and invariant hermitian
structures $(J, g )$ on $T^*G$, where $g$ is the standard neutral
metric on $T^*G$. In \S\ref{no-GCS} we prove Theorem~\ref{main}
which asserts that
\emph{a four dimensional solvable Lie group $G$ has neither left
invariant symplectic nor complex structures if and only if $G$ does not admit
generalized complex structures.} In the proof, we use
the classification
of 4--dimensional solvable Lie groups with
left invariant complex (resp. symplectic) structures
carried out in~\cite{SJ}
and~\cite{O1} (resp.~\cite{MR}; see also~\cite{O}).

On the other hand, in \S\ref{sec-5} we distinguish the solvable
Lie groups of dimension $4$ admitting a non-extremal left
invariant generalized complex structure (\S\ref{subsec-5}) and the
Lie groups carrying a left invariant complex or symplectic
structure but without a non-extremal left invariant generalized
complex structure (\S\ref{sec6}).

Finally, in \S\ref{six-dim} we show that Theorem~\ref{main} does not work
in dimension $6$. In fact, we construct
an example of a six dimensional (non-nilpotent) solvable Lie group admitting neither
left invariant symplectic  nor
complex structures but having non-extremal generalized complex structures.

\section{Hermitian structures on cotangent Lie groups}\label{prelim}
A \emph{left invariant complex structure} on a real Lie group $G$ is
a complex structure on the underlying manifold  such that left
multiplication  by elements of the group are holomorphic.
Equivalently, there exists an endomorphism $J$ of $\mathfrak g$, the
Lie algebra of $G$, such that: \\ $J^2 = -\id$ and $N_J \equiv 0$,
where
\begin{equation}\label{nijen}
N_J(x,y)=   [x,y] + J[Jx,y] + J[x,Jy] -[Jx,Jy], \qquad\forall\,x,y \in\mathfrak g .
\end{equation}
The condition $N_J \equiv 0$ is called the \emph{integrability
condition} of $J$.

The action of $G$ on itself given by left multiplication $L_g\colon
G \longrightarrow G \,$ can be lifted to an action of $G$ on $TG$
given by $dL_g\colon TG\longrightarrow TG$.  Thus, a left invariant
complex structure is an equivariant endomorphism of $TG$ with
respect to the lifted action of $G$ given by left multiplication.
Similarly, a left invariant symplectic structure on $G$ is an
equivariant isomorphism $\omega\colon TG\longrightarrow T^*G$ where
the action of $G$ on $T^*G$ is $L_{g^-1}^* \colon T^*G
\longrightarrow T^*G$.

A \emph{left invariant hermitian structure} on $G$ is a pair $(J,g)$ of a
left invariant complex structure $J$ together with a left invariant
hermitian metric $g$ (not necessarily positive definite). If $J$ denotes
the corresponding endomorphism on $\mathfrak g$ and $\langle \cdot \, ,
\cdot \rangle$ the non degenerate symmetric bilinear form on $\mathfrak
g$ induced by $g$,   we say that
$(J, \langle \cdot \, , \cdot \rangle)$ is a hermitian structure on $\mathfrak g$. A
non degenerate symmetric bilinear form $\langle \cdot\, ,\cdot \rangle$ on
$\mathfrak g$ is said to be ad-invariant when it satisfies:
\begin{equation}\label{bi-inv}
\langle [x,y],z\rangle + \langle y,[x,z] \rangle = 0 \qquad \text{ for any } x, y, z
\in \mathfrak g .
\end{equation}
If $G$ is a Lie group with Lie algebra $\mathfrak g$ and $g$ is a bi-invariant metric
on $G$, that is, $g$ is both left and right invariant, then the bilinear form
$\langle\cdot \, ,\cdot\rangle$ on $\mathfrak g$ induced by $g$ is ad-invariant.

Let $\mathfrak g^{\mathbb C} =\mathfrak g \otimes_{\mathbb R} \mathbb C$  be the
complexification of the real Lie algebra $\mathfrak  g$ and let $\sigma$ denote the
conjugation in $\mathfrak  g^{\mathbb C}$ with respect to the real form $\mathfrak
g$, that is, $\sigma (x + i y) = x - i y$, $x, y\in \mathfrak g$. Starting with a
hermitian structure $(J, \langle \cdot\,,\cdot \rangle)$ on $\mathfrak g$, let
$J^{\mathbb C}$ (resp. $\langle\cdot\, ,\cdot \rangle^{\mathbb C}$) denote the
complex linear (resp. complex bilinear) extension of $J$ (resp. $\langle\cdot\,,
\cdot\rangle$) to $\mathfrak  g^{\mathbb C}$. We obtain a splitting
\[
\mathfrak g^{\mathbb C} = \mathfrak q \oplus\sigma (\mathfrak q),
\]
where $\mathfrak q$, the $i$-eigenspace of $J^{\mathbb C}$, is a complex subalgebra
of $\mathfrak g^{\mathbb C}$ which is maximal isotropic with respect to $\langle
\cdot\,,\cdot\rangle^{\mathbb C}$.

We prove the above statement in the  following proposition, where it is shown that,
conversely, if $B$ is a  symmetric bilinear form on $\mathfrak g^{\mathbb C}$
satisfying certain conditions, then any splitting  $\mathfrak g^{\mathbb C} =
\mathfrak q\oplus \sigma (\mathfrak q)$, where $\mathfrak q$ is a maximal
$B$-isotropic complex subalgebra of $\mathfrak g^{\mathbb C}$, gives rise to a
hermitian structure $(J,\langle\cdot\,,\cdot\rangle)$ on $\mathfrak g$ such that
the $i$-eigenspace of $J^{\mathbb C}$ in $\mathfrak g^{\mathbb C}$ is $\mathfrak q$
and $\langle x, y \rangle = B(x,y)$ for $x,y \in \mathfrak  g$.

\begin{prop}\label{prop2}
Let $G$ be a Lie group with Lie algebra $\mathfrak g$ and denote by $\mathfrak
g^{\mathbb C}$ the complexification of $\mathfrak g$. There is a one-to-one
correspondence between left invariant hermitian structures $(J,g)$ on  $G$  and
pairs $(\mathfrak q,B)$, where  $B$ is a symmetric bilinear form  on  $\mathfrak
g^{\mathbb C}$ and $\mathfrak q$ is a maximal $B$-isotropic complex subalgebra of
$\mathfrak g^{\mathbb C}$ satisfying the following conditions:
 \begin{eqnarray}
\mathfrak  g^{\mathbb C} &=& \mathfrak q \oplus \sigma (\mathfrak q),\label{def2}\\
\label{real} B(\sigma z,\sigma w)&=& \overline{B(z,w)}, \qquad z,w\in\mathfrak
g^{\mathbb C},
\end{eqnarray}
where $\overline{\alpha }$ denotes the complex  conjugate of $\alpha\in\mathbb C$
and $\sigma$ is the conjugation in  $\mathfrak g^{\mathbb C}$ with respect to
$\mathfrak  g$.
\end{prop}
\begin{proof}   Given a left invariant hermitian structure $(J,g)$ on $G$, let $(J,
\langle \cdot \, , \cdot \rangle)$ be the corresponding  hermitian structure on
$\mathfrak  g$. $\mathfrak g^{\mathbb C}$ decomposes into a direct sum of subspaces
$\mathfrak  g^{\mathbb C}=\mathfrak g^{1,0}\oplus\mathfrak g^{0,1} $, the eigenspaces
of $J^{\mathbb C}$ of eigenvalue $i$ and $-i$, respectively. It follows that
\[
\mathfrak g^{1,0}=  \{ x-iJx : x \in \mathfrak  g \}, \qquad \qquad
\mathfrak  g^{0,1}=  \{ x+iJx : x \in \mathfrak  g\},
\]
hence, $\mathfrak  g^{0,1}=\sigma (\mathfrak g^{1,0})$.
Equation~$N_J \equiv 0$ is equivalent to the fact that these
subspaces are subalgebras. Moreover, using that $J$ is orthogonal,
it is easy to check that both $\mathfrak g^{1,0}$ and $\mathfrak
g^{0,1}$  are isotropic with respect to $\langle \cdot \, ,\cdot
\rangle^{\mathbb C}$. Since $\langle \cdot \, , \cdot \rangle$ is
non degenerate, these subalgebras are maximal isotropic. Hence,
$(\mathfrak g^{0,1},\langle \cdot \,, \cdot \rangle^{\mathbb C})$
satisfies the required conditions. Note that equation \eqref{real}
holds if and only if $B$ takes real values on $\mathfrak g$, and
$\langle \cdot \, , \cdot \rangle ^{\mathbb C}$ clearly satisfies
this property.

Conversely, given a  pair $(\mathfrak q, B)$ as in the statement, we wish to show
that it gives rise to a hermitian structure $(J,\langle \cdot \, , \cdot \rangle )$
on $\mathfrak  g$. Let $J$ be the almost complex structure defined on $\mathfrak
g^{\mathbb C}$ by
\[
Jz = i z,  \qquad \quad J\circ\sigma (z) = -i \, \sigma (z), \qquad z \in \mathfrak q
.
\]
Since $J \circ \sigma = \sigma\circ  J$, then $J$  leaves
$\mathfrak g$ stable. The fact that $\mathfrak q$ is a subalgebra
implies that $J$ satisfies $N_J \equiv 0$. Since equation
\eqref{real} holds, $B$ takes real values on $\mathfrak g$. Let
$\langle \cdot \, , \cdot\rangle$ be the restriction of $B$ to
$\mathfrak g$. It follows from \eqref{def2} and the fact that
$\mathfrak q$ is $B$-isotropic  that $J$ is orthogonal with
respect to $\langle \cdot \, , \cdot \rangle$. Since  $\mathfrak
q$ is maximal isotropic then $\langle \cdot \, , \cdot \rangle$ is
non degenerate, that is,  $(J, \langle \cdot \, ,\cdot \rangle )$
is a hermitian structure on $\mathfrak g$. Therefore, it induces,
by left translations, a left invariant hermitian structure on $G$
and the proposition follows.
\end{proof}


We will be studying a special class  of left invariant hermitian structures.
The Lie groups that come into the picture are the cotangent bundles of  Lie groups
with a standard bi-invariant metric.

Let $\mathfrak g$ be a Lie algebra and $\mathfrak v$ a $\mathfrak g$-module, that is,
there exists a Lie algebra homomorphism $\rho\colon\mathfrak g
\longrightarrow \mathfrak{gl}(\mathfrak v)$. Let  $\mathfrak g\ltimes_{\rho}
\mathfrak v$ denote the  semidirect product of $\mathfrak g$ by $\mathfrak v$, where
we look upon $\mathfrak v$ as an abelian Lie algebra. The bracket on  $\mathfrak g
\ltimes_{\rho}\mathfrak v$  is given as follows:
\begin{equation} \label{bracket2}
[(x,u),(y,v)]=([x,y],\rho(x)v -\rho(y)u) \hspace{1cm} \text{ for }x, y\in
\mathfrak g, \; u,v \in \mathfrak v.
\end{equation}
Complex structures on Lie algebras of the above type were studied in \cite{BD}. In
the present article we will restrict our attention to the particular case when
$\mathfrak v=\mathfrak g^*$ and $\rho= \mbox{ad}^*$ is the coadjoint  representation:
\[
\mbox{ad}^* (x) (\alpha)= -\alpha \circ \mbox{ad}(x), \qquad \alpha \in\mathfrak g^*,
\; x \in \mathfrak g .
\]
We will denote $\mathfrak g \ltimes_{\mbox{ad}^*}\mathfrak g^*$ by $(T^*\mathfrak
g,\mbox{ad}^*)$, the Lie bracket being given by \eqref{bracket1}.
The cotangent algebra $(T^*\mathfrak g , \text{ad}^*)$ has a standard non degenerate
symmetric ad-invariant  bilinear form $\langle \cdot \, ,\cdot \rangle$  (see
\eqref{st-bilform}). We notice that the subalgebra $\mathfrak g$ and the ideal
$\mathfrak g^*$ are maximal isotropic in $(T^*\mathfrak g,\langle \cdot \,,\cdot
\rangle)$.

Left invariant hermitian structures on the cotangent Lie group $T^*G$ are given by
endomorphisms $J$ of $T^*\mathfrak g$ whose matrix form with respect to the
decomposition $\mathfrak g \oplus \mathfrak g^*$ is
\[
J = \begin{pmatrix}J_1& J_2\\  J_3&J_4\end{pmatrix},
\]
and satisfy
\begin{equation}\label{formadeJ-10} \begin{split}
\mbox{(i)} \; &  J_4 = -  J_1^*,\quad  J_2=-J_2^*,\quad  J_3=- J_3^*,\\
\mbox{(ii)} \; & J_1^2+ J_2 J_3 =-\id,\quad J_1J_2=-(J_1J_2)^*,\quad
J_3J_1=-(J_3J_1)^*,\\
\mbox{(iii)} \; & J \text{ is integrable.} \end{split}
\end{equation}

\begin{ex} \label{type-n}
Let $\jmath$ be a complex structure on $\mathfrak g$, $\dim\mathfrak g =2n$, and
define $J_{\jmath}$ on $T^*\mathfrak g$ by
\begin{equation}\label{tipo-n}
J_{\jmath}(x,\alpha)=(\jmath(x),- \jmath^*(\alpha )), \qquad x\in \mathfrak g , \;
\alpha \in \mathfrak g^*,
\end{equation}
where $\jmath^*$ is the adjoint of $\jmath$, that is, $\jmath^*(\alpha) =\alpha
\circ\jmath$. It follows that $J_{ \jmath}$ is orthogonal with respect to the
standard bilinear form $\langle \cdot \, ,\cdot \rangle$ on $T^*\mathfrak g$.
Moreover, it was shown in \cite{BD} (Proposition 3.2) that the integrability
of $\jmath$ implies that $J_{ \jmath}\,$ is a complex structure on $(T^*\mathfrak g,
\text{ad}^*)$. Therefore, $(J_{\jmath},\langle \cdot \, ,\cdot \rangle)$ is a
hermitian structure on $(T^*\mathfrak g, \text{ad}^*)$.
\end{ex}

\begin{ex}\label{type-0}
Let $\omega\colon\mathfrak g\longrightarrow\mathfrak g^*$ be a linear isomorphism and
define
\begin{equation}\label{tipo-0}
 J_\omega(x,\alpha)= (-\omega^{-1}(\alpha),\omega(x)),
\end{equation}
(compare with \S4 in \cite{BD}). It follows that $J_{\omega}$ is orthogonal with
respect to the standard bilinear form on $T^*\mathfrak g$ if and only if $\omega$ is
skew-symmetric. The integrability of $J_{\omega}$ is equivalent to the following
condition
\begin{equation} \label{CYBE}
\omega ([x,y])= \omega (x)\circ\text{ad}(y) -\omega (y) \circ \text{ad} (x) .
\end{equation}
Therefore, if $\omega$ satisfies \eqref{CYBE} $J_{\omega}$ defines a hermitian
structure on $(T^*\mathfrak g,\text{ad}^*)$. We observe that in this case, $\omega$
is  a symplectic structure  on $\mathfrak  g$.
\end{ex}


\section{Left invariant generalized complex structures on Lie groups}\label{cotg}
We recall that a generalized complex structure on a manifold $M$ is an endomorphism
$\mathcal J$ of $TM \oplus T^*M$ satisfying $\mathcal J^2=-\id$ which is orthogonal
with respect to the standard inner product $\langle \cdot \, ,\cdot\rangle$ on $TM
\oplus T^*M$  defined in \eqref{st-bilform} and such that the $i$ eigenbundle of
$\mathcal J$ in $(TM \oplus T^*M)\otimes \mathbb C$ is involutive with respect to
the Courant bracket. This bracket is defined as follows:
\begin{equation}\label{courant}
\bigr[(X,\xi), (Y,\eta)\bigl]= \bigr([X,Y],{\mathcal L}_X\eta- {\mathcal L}_Y\xi- \frac
12 d(i_X\eta -i_Y \xi)\bigl),
\end{equation} where $(X,\xi)$, $(Y,\eta)$ are smooth sections of $TM\oplus T^*M$.

When $M$  is a Lie group $G$, consider the left action of $G$ on $TG\oplus T^*G$
induced by left multiplication of $G$ on itself, that is,
\begin{equation} \label{action}
\begin{split}
\lambda\colon G\times(TG\oplus T^*G) &\to TG \oplus T^*G,\\
\left(g ,(x,\alpha)\right)\hspace{.5cm}& \mapsto \left((dL_{g})_h\, x,
(L_{g^{-1}}^*)_{gh} \, \alpha \right) , \quad x\in T_hG, \;\; \alpha \in T_h^*G,
\;\; g, h\in G
\end{split}
\end{equation}
where
\[
(L_{g^{-1}}^*)_{gh} \, \alpha \, (y)= \alpha\left((dL_{g^{-1}})_{gh} \, y
\right), \qquad\forall \, y \in  T_{gh}G .
\]
A generalized complex structure $\mathcal J$ on $G$ is said to be left invariant (or
$G$-invariant) if
\[
\mathcal J\colon TG \oplus T^*G \longrightarrow TG\oplus T^*G
\]
is equivariant with respect to the induced left action of $G$ on $TG\oplus T^*G$
given in \eqref{action}. It follows that, for any $g\in G$, the following diagram is
commutative:
\[
\begin{CD} T_gG\oplus T_g^*G @>{\mathcal J_g}>>  T_gG \oplus T_g^*G\\
@V{\lambda _{g^{-1}}}VV              @VV{\lambda_{g^{-1}}}V\\
\mathfrak g\oplus \mathfrak g^* @>>{\mathcal J_e}>\mathfrak g
\oplus \mathfrak g^*
\end{CD},
\]
where
\[
\lambda_{g^{-1}}(x,\alpha) =\lambda\left(g^{-1},(x,\alpha)\right), \qquad x\in T_gG,
\;\; \alpha \in T_g^*G .
\]
In other words, $\mathcal J$ is left invariant if and only if, for any $g \in G$,
$\mathcal J_g$ is given in terms of $\mathcal J_e$ as follows:
\begin{equation}\label{inv-gencx}
\mathcal J_g= \lambda_g\circ\mathcal J_e \circ\lambda_{g^{-1}}=
\begin{pmatrix}(dL_{g})_e & \\ & (L_{g^{-1}}^*)_g\end{pmatrix}
\circ\mathcal J_e\circ\begin{pmatrix} (dL_{g^{-1}})_g & \\& (L_{g}^*)_e
\end{pmatrix}.
\end{equation}
If we identify the space of left invariant sections of $\, TG\oplus
T^*G\, $ with $\, \mathfrak g\oplus\mathfrak g^*$, then the
restriction of the Courant bracket \eqref{courant} to $\mathfrak
g\oplus\mathfrak g^*$ is precisely the Lie bracket \eqref{bracket1}
on the cotangent algebra $(T^*\mathfrak g ,\text{ad}^*)$. Therefore,
the Courant integrability condition of a left invariant generalized
complex structure $\mathcal J$ on $G$ is equivalent to the
integrability of $\mathcal J_e$ on the cotangent algebra
$(T^*\mathfrak g,\text{ad}^*)$. Moreover, since $\lambda_g$, $g \in
G$, are isometries of the standard bilinear form $\langle \cdot \, ,
\cdot\rangle$ on $TG\oplus T^*G$, it follows that $\mathcal J$ is
orthogonal with respect to $\langle\cdot \, , \cdot \rangle$ if and
only if $\mathcal J_e$ is compatible with $\langle\cdot \, ,\cdot
\rangle$. Therefore, if $\mathcal J$ is a left invariant generalized
complex structure on $G$, $(\mathcal J_e,\langle\cdot\,
,\cdot\rangle)$ is a hermitian structure on $T^*\mathfrak g$.
Conversely, given a hermitian  structure $(J,\langle\cdot\,
,\cdot\rangle)$ on $(T^*\mathfrak g,\text{ad}^*)$, where
$\langle\cdot\, ,\cdot\rangle$ is the standard neutral metric on
$T^*\mathfrak g$, it can be extended, by means of \eqref{inv-gencx},
to a left invariant generalized complex structure $\mathcal J$ on
$G$ such that $\mathcal J _e =J$.

The preceding arguments  yield the following result:

\begin{prop}\label{corresp}
There is a one-to-one correspondence between left invariant generalized complex
structures on $G$ and invariant hermitian structures $(J, g )$ on $T^*G$, where $g$
is the standard neutral metric on $T^*G$.
\end{prop}

When a Lie group $G$ has a left invariant complex  or symplectic structure, then any
of these structures induces a natural left invariant generalized complex structure
on $G$, as shown in Examples \ref{type-n} and \ref{type-0}.

In view of Proposition \ref{corresp}, a hermitian structure on $(T^*\mathfrak g
,\text{ad}^*)$ with respect to the standard bilinear form  will be called a
generalized complex structure on $\mathfrak g$ and denoted by $(\mathcal J,\langle
\cdot\, ,\cdot\rangle)$. When $\mathcal J$ satisfies only conditions (i) and (ii) in
\eqref{formadeJ-10}, it will be called an almost generalized complex structure. Note
that if $T^*\mathfrak g$ is a generalized complex vector space, $\dim \mathfrak g
=2n$ (see \cite{Gualt,BBB}).

\begin{remark}
It was proved in \cite{FGG}  that every four dimensional nilpotent
Lie group has either left invariant complex or symplectic structures
(maybe both; see also \cite{Mag} for the classification of these
groups). Hence, such a Lie group has a left invariant generalized
complex structure. In \cite{CG} (see also \cite{Cav}) it was shown
that every six dimensional nilpotent Lie group admits a left
invariant generalized complex structure.  In other words, the
cotangent algebra $(T^*\mathfrak g ,\text{ad}^*)$ of any four or six
dimensional nilpotent Lie algebra $\mathfrak g$ admits a hermitian
structure $(J,\langle\cdot\, ,\cdot\rangle)$,  where $\langle\cdot\,
,\cdot\rangle$ is  the standard bilinear form  on $T^*\mathfrak g$.
\end{remark}

It was proved in \cite{Ale} that when $(J,g)$ is a left invariant
hermitian structure on a Lie group $H$ such that $g$ is bi-invariant
then both $H$ and $H^*$ are Poisson Lie groups. Moreover, since
$J$ is a complex structure, $H^*$ is a complex Lie group (see
\cite{ABO}).  The Poisson structure on $H$ is given by $\Pi(h)=
(dR_h)_eJ - (dL_h)_eJ$, $h\in H$,  where $J$ is viewed as an element
of $\mathfrak h\wedge\mathfrak h$ by identifying the Lie algebra
$\mathfrak h$ of $H$ with its dual $\mathfrak h^*$ via the metric
$g$.  As a corollary of this result and Proposition~\ref{corresp},
we therefore obtain:

\begin{cor}\label{Poisson}
If $G$ is a Lie group with a left invariant generalized complex structure, then
$T^*G$ and $(T^*G)^*$ are Poisson Lie groups such that $(T^*G)^*$ is a complex Lie
group.
\end{cor}

We end this section by determining the generalized complex structures on the two
dimensional non-abelian Lie algebra $\mathfrak g = \mathfrak{aff}(\mathbb R)$.

\begin{ex}\label{cot-aff(R)}
\emph{Generalized complex structures on $\mathfrak{aff}(\mathbb R)$.}
Let $\mathfrak g =\mathfrak{aff}(\mathbb R)$ be the two dimensional non-abelian Lie
algebra and $T^*\mathfrak{aff}(\mathbb R)$ the corresponding cotangent Lie algebra.
Let $\{e_0,e_1\}$ be a basis of $\mathfrak g$ such that $[e_0 ,e_1] =e_1$, and
$\{\alpha^0,\alpha ^1\}$ the dual basis of $\mathfrak g^*$. Set
\[
X_i=(e_{i-1},0), \qquad X_{i+2}=(0,\alpha^{i-1}), \qquad i=1,2,
\]
then:
\[
[X_1 ,X_2] =X_2, \qquad [X_1 ,X_4] =-X_4, \qquad [X_2 ,X_4] =X_3 .
\]
A generalized complex structure $\mathcal J$ on $\mathfrak{aff}(\mathbb R)$ takes the
following form in the ordered basis $\{X_1, \ldots , X_4\}$:
\[
\mathcal J =
\begin{pmatrix}
a_{11} & a_{12}& 0& a_{14}\\ a_{21} & a_{22}& -a_{14}& 0\\
0& -a_{41}& -a_{11}&-a_{21}\\a_{41}&0&-a_{12}&-a_{22}
\end{pmatrix},
\]
with $\mathcal J ^2= -\id$ and $N_{\mathcal J}\equiv 0$.

In case $a_{14}\neq 0$, the condition $\mathcal J^2= -\id$ implies $a_{41}
\neq 0$, $a_{11}^2 +a_{14}\, a_{41}=-1$, $a_{11}=a_{22}$ and $a_{12}=0=a_{21}$.
Hence,
\begin{equation}\label{aff(R)-0}
\mathcal J = \begin{pmatrix}
a_{11} & 0& 0& a_{14}\\0 & a_{11}& -a_{14}& 0\\
0& -a_{41}& -a_{11}&0\\a_{41}&0&0&-a_{11}
\end{pmatrix},
\qquad a_{14}\, a_{41}\neq 0, \quad a_{11}^2+a_{14}a_{41}=-1.
\end{equation}
It follows that $\mathcal J$ as above satisfies $N_{\mathcal
J}\equiv 0$. In particular, if $a_{11}=0$, $\mathcal J$ arises from
a symplectic structure on $\mathfrak{aff}(\mathbb R)$ as in
Example~\ref{type-0}, but for $a_{11}\neq 0$ $\mathcal J$  is not
induced by a symplectic or complex structure on
$\mathfrak{aff}(\mathbb R)$. However, since $\mathcal J$ is of type
0 (see the paragraph next to \eqref{diag} in \S\ref{no-GCS}) it
follows from Theorem~\ref{simple} that it is equivalent to a
symplectic structure via a $B$-field transformation.

In case $a_{14}=0$, the condition $\mathcal J ^2= -\id$ implies
\[ a_{41}=0,\quad a_{11}^2 +a_{12}a_{21}=-1,\quad a_{12}\, a_{21} \neq
0, \quad a_{11}=-a_{22}. \] Therefore,
\begin{equation}\label{aff(R)-1}
\mathcal J = \begin{pmatrix}
a_{11} & a_{12}& 0&0 \\a_{21} & -a_{11}& 0& 0\\
0& 0& -a_{11}&-a_{21}\\0&0&-a_{12}&a_{11}
\end{pmatrix},
\qquad a_{12}\, a_{21}\neq 0, \quad a_{11}^2+a_{12}\, a_{21}=-1 ,
\end{equation}
and $\mathcal J$ satisfies $N_{\mathcal J}\equiv 0$. Note that every generalized
complex structure in this family arises from a complex structure on
$\mathfrak{aff}(\mathbb R)$ as in Example~\ref{type-n}.

We observe that $T^*\mathfrak{aff}(\mathbb R)$ is isomorphic to the Lie algebra
$\mathfrak d_4$ (see \cite{ABDO}). This is the unique four dimensional solvable Lie
algebra admitting a structure of a Manin triple. The above calculations together
with Corollary \ref{Poisson} imply that the Lie group $\mathcal D_4$ with Lie
algebra $T^* \mathfrak{aff}(\mathbb R)$ is a Poisson Lie group such that the
Poisson Lie group ${\mathcal D_4}^*$ is a complex Lie group.

Fix two generalized complex structures $\mathcal J_1$, $\mathcal J_2$ on $\mathfrak
g=\mathfrak{aff}(\mathbb R)$ as follows:
\[
\mathcal J_1=\begin{pmatrix}
a & 0& 0& b\\0 & a& -b& 0\\0& -c& -a&0\\c&0&0&-a
\end{pmatrix}, \ a^2+bc=-1,
\qquad
\mathcal
J_2=\begin{pmatrix}
x & y& 0&0 \\z & -x& 0& 0\\0& 0& -x&-z\\0&0&-y&x
\end{pmatrix}, \ x^2+yz=-1,
\]
and consider $G=-\mathcal J_1\mathcal J_2$. Observe that $\mathcal
J_1$ and $\mathcal J_2$ commute, therefore $G^2=\id$. It follows
that $G$ defines a positive definite metric on $\mathfrak
g\oplus\mathfrak g^*$ if and only if $cz<0 $. Therefore, when this
 condition is satisfied, we obtain generalized K\"ahler
structures on $\mathfrak{aff}(\mathbb R)$ (see \cite{Cav,Gualt}).
\end{ex}

\section{Solvable Lie groups without
generalized complex structures}\label{no-GCS}

In this section we prove that a four dimensional (non-nilpotent)
solvable Lie group has no left invariant generalized complex structures
if and only if it admits neither left invariant symplectic nor left
invariant complex structures.




We start by fixing some notation. Let $\{\alpha^i\}_{i=0}^3$ be the basis of $\mathfrak
g^*$ dual to the basis $\{e_i\}_{i=0}^3$ of $\mathfrak g$. Define the basis
$\left\{X_i\right\}_{i=1}^8$ of $T^*\mathfrak g$ by
\begin{equation}\label{base}
X_i= (e_{i-1},0)\qquad \mbox{and} \qquad
X_{i+4}=(0,\alpha^{i-1}),\qquad 1\leq i \leq 4.
\end{equation}
Let $\mathcal J$ be a linear endomorphism of $T^*\mathfrak g$ whose matrix form is
\begin{equation}\label{jotas4}
\mathcal J = \begin{pmatrix}
\mathcal J_1&\mathcal J_2\\
\mathcal J_3&\mathcal J_4
\end{pmatrix},
\end{equation}
with respect to the basis  $\left\{X_i\right\}_{i=1}^8$ of $T^*\mathfrak g$ defined
by~\eqref{base}. If  $\mathcal J$ is orthogonal with respect to the standard bilinear
form on $T^*\mathfrak g$  then  the matrix of $\mathcal J$ is of the form

\begin{equation}\label{jota88}
\mathcal J=\begin{pmatrix}
a_{11} &a_{12} &a_{13} &a_{14}  &0      &a_{16} &a_{17} &a_{18}  \\
a_{21} &a_{22} &a_{23} &a_{24}  &-a_{16}&0      &a_{27} &a_{28}  \\
a_{31} &a_{32} &a_{33} &a_{34}  &-a_{17}&-a_{27}&0      &a_{38}  \\
a_{41} &a_{42} &a_{43} &a_{44}  &-a_{18}&-a_{28}&-a_{38}&0       \\
0      &-a_{61}&-a_{71}&-a_{81} &-a_{11}&-a_{21}&-a_{31}&-a_{41} \\
a_{61}&0       &-a_{72}&-a_{82} &-a_{12}&-a_{22}&-a_{32}&-a_{42} \\
a_{71}&a_{72}  &0      &-a_{83} &-a_{13}&-a_{23}&-a_{33}&-a_{43} \\
a_{81}&a_{82}  &a_{83} &0       &-a_{14}&-a_{24}&-a_{34}&-a_{44}
\end{pmatrix}.
\end{equation}
Moreover, taking into account \eqref{formadeJ-10}, if $\mathcal J^2=-\id$ then the matrix
$\mathcal J$ has the following property:
\begin{equation}\label{diag}
\text{ for every }
1\leq i\leq 4  \text{ there exists } j\neq i \text{ such that } a_{ij}\not=0.
\end{equation}

\medskip
We will say that $\mathcal J$ is of {\em complex type} if $\, \mathcal J_2
=\mathcal J_3=0\, $, $\mathcal J$ is of {\em symplectic type} if $\, \mathcal J_1
=\mathcal J_4=0\,$, and  $\mathcal J$ is said to be of {\em type~$k$} when
$\rank(\mathcal J_2) =2(n-k)$, where $\dim \mathfrak g= 2n$ (compare with
\cite{Gualt}). Observe that if $\mathcal J$ is of complex (resp.
symplectic)  type then it is of type~$2$ (resp.~$0$).
\medskip

We recall a theorem from~\cite{Gualt,CG}

\begin{thm}\label{simple}\emph{(\cite{CG}, Theorem 1.1; \cite{Gualt}, Theorem 4.35)}
Any regular point of type $k$ in a generalized complex
$2n$-manifold has a neighbourhood which is equivalent, via a
diffeomorphism and a $B$-field transformation, to the product of an
open set in $\mathbb C^k$ with an open set in the standard
symplectic space $\mathbb R^{2n-2k}$.
\end{thm}

The previous theorem implies that a $2n$-dimensional Lie
algebra admits a generalized complex structure of type 0 (resp. of
type $n$)  if and only if it has a symplectic structure (resp.,
a complex structure).



In order to prove the main result of this section, we
recall the definition of the four dimensional solvable Lie algebras admitting
neither symplectic nor complex structures (see \cite{Sam,Sa,MR}). They are
\begin{equation}\label{cincofam}
\begin{split}
\mathbb R\times \mathfrak r_3:\ &[e_1, e_2]=e_2, \;\;  [e_1,e_3]=e_2+ e_3\\
\mathbb R \times \mathfrak r_{3, \lambda}:\ & [e_1, e_2]=e_2, \;\;
[e_1,e_3]=\lambda e_3,\;\;|\lambda|< 1 , \;\; \lambda \neq  0;\\
\mathfrak r_{4}:\ &[e_0,e_1]=e_1, \;\; [e_0, e_2]=e_1+e_2, \;\;
[e_0,e_3]=e_2+e_3;\\
\mathfrak r_{4, \lambda}:\ & [e_0, e_1]=e_1, \;\;
[e_0,e_2]=\lambda e_2,\;\; [e_0, e_3]= e_2+\lambda e_3,\;\;\lambda\in\mathbb R,\;\;
\lambda \neq -1,\, 0,\, 1;\\
\mathfrak r_{4,\mu,\lambda}:\ &[e_0,e_1]=e_1, \; \;
[e_0, e_2]= \mu e_2, \;\; [e_0,e_3]= \lambda e_3,
\;\;
 -1 < \mu <\lambda  <1,\\   & \hspace{8cm}\mu \lambda \neq 0,\;\;\mu +\lambda \neq0.
\end{split}
\end{equation}

Next, we show that every Lie algebra $\mathfrak h$ included in
the list \eqref{cincofam} has no left invariant generalized complex structures by
analyzing each case. To this end, we will prove that any almost complex  structure
$\mathcal J$ on $T^*(\mathfrak h)$ does not satisfy the integrability condition.
This condition is equivalent to the vanishing of the 256 coefficients
$N_{ij}^k$ defined by
\[
N_{\mathcal J}(X_i,X_j)=\sum_{k=1}^8 N_{ij}^k X_k, \qquad 1\leq i< j\leq 8 ,
\]
where $N_{\mathcal J}$  is the Nijenhuis tensor of $\mathcal J$ (see \eqref{nijen}).


\begin{prop} \label{Lie1}
The Lie algebra $\mathbb R\times\mathfrak r_3$ does not admit generalized complex
structures.
\end{prop}

\begin{proof}
We consider  the basis $\{X_i\}$ for $T^*(\mathbb
R\times\mathfrak r_3)$ defined by \eqref{base} .
Taking into account \eqref{courant} and the definition of the
Lie algebra $\mathbb R\times\mathfrak r_3$ given in \eqref{cincofam}, we see
that the only non-zero Lie brackets on $T^*(\mathbb R\times\mathfrak r_3)$ are
\[
\begin{array}{lll}
[X_2,X_3]=X_3,\quad &[X_2,X_4] =X_3+X_4,&[X_2,X_7]=-X_7-X_8,\\[5pt]
[X_2,X_8]=-X_8, &[X_3,X_7]=X_6,&[X_4,X_7]=X_6=[X_4,X_8].
\end{array}
\]

Suppose that $\mathbb R\times\mathfrak r_3$ has a generalized complex structure,
i.e., $T^*(\mathbb R\times\mathfrak r_3)$ has a hermitian structure $(\mathcal J,
\langle\cdot\, ,\cdot\rangle)$. Since all the coefficients $N_{ij}^k$ of the
Nijenhuis tensor $N_{\mathcal J}$ of $\mathcal J$ must be zero, we have $0=N_{78}^2=a_{28}{}^2$
and $0=N_{46}^7=a_{23}{}^2$, and so $a_{28}=a_{23}=0$.

Let us consider the equation
\begin{equation}\label{r344}
0=N_{48}^6=1+a_{44}{}^2+a_{43}(a_{22}+a_{44}+a_{34})+a_{24}a_{42}-a_{38}a_{83}.
\end{equation}

Now, \eqref{r344} and the equations
\[
0=N_{67}^8=2 a_{24} a_{27},\qquad 0=N_{68}^6=-a_{43} a_{27},\qquad
0=N_{37}^8=-a_{43} a_{24}-2 a_{83} a_{27},
\]
imply that $a_{27}=0$. Moreover, because
\[
0=N_{78}^1=-2 a_{16} a_{38},\qquad0=N_{37}^1=a_{16} a_{43},\qquad
0=N_{46}^1=a_{16} a_{24},\]
we have $a_{16}=0$ using \eqref{r344}; and because
\begin{align*}
&0=N_{37}^8=- a_{24} a_{43},\qquad 0=N_{48}^3=-2a_{24} a_{38},\quad\mbox{and}\\
&0=N_{27}^7=-1-a_{33}{}^2+a_{43}(a_{22}-a_{33}-a_{34})-a_{38}a_{83},
\end{align*}
we see that $a_{24}=0$.  Finally,
\[
0=N_{78}^5=2a_{21} a_{38},\qquad 0=N_{18}^8=-a_{21} a_{43}
\]
imply $a_{21}=0$ using \eqref{r344}. So, in the
matrix \eqref{jota88} of $\mathcal J$, the unique  non-zero entry in the 2nd row is
$a_{22}$, which is not possible  by \eqref{diag}. This shows that
$\mathcal J$ cannot be integrable.
\end{proof}


\begin{prop} \label{Lie2}
For $\lambda\not=0,\pm1$, the Lie algebra $\mathbb R\times\mathfrak r_{3,\lambda}$
has no generalized complex structure.
\end{prop}

\begin{proof}
Using the basis $\{X_i\}$ for $T^*(\mathbb
R\times\mathfrak r_{3,\lambda})$, given by  ~\eqref{base}, and the definition of
the Lie algebra $\mathbb R\times\mathfrak r_{3,\lambda}$ stated in~\eqref{cincofam},
the only non-zero Lie brackets on $T^*(\mathbb R\times\mathfrak r_{3,\lambda})$ are
\[
\begin{array}{lll}
[X_2,X_3]=X_3,&[X_2,X_4]=\lambda X_4,&[X_2,X_7]=-X_7\\[5pt]
[X_2,X_8]=-\lambda X_8,&[X_3,X_7]=X_6,&[X_4,X_8]=\lambda X_6.
\end{array}
\]

Let $(\mathcal J,\langle\cdot\, ,\cdot\rangle)$ be a hermitian structure on
$T^*(\mathbb R\times\mathfrak r_{3,\lambda})$. The integrability condition
of $\mathcal J$ implies that all the coefficients $N_{ij}^k$
of the Nijenhuis tensor of $\mathcal J$ are zero. In particular,
\begin{equation}\label{r3l44}
0=N_{48}^6=\lambda(1+a_{44}{}^2+a_{24}a_{42}+a_{28}a_{82})+ a_{34}a_{43}-a_{38}a_{83}
.
\end{equation}
Now from \eqref{r3l44} and the equations
\begin{align*}
&0=N_{68}^1= -\lambda a_{16} a_{28},  \qquad 0=N_{45}^2=\lambda a_{16}a_{24}, \qquad
0=N_{57}^2= a_{16} a_{27},\\
&0=N_{57}^8= \lambda a_{14}a_{27}+a_{17}a_{24}+(\lambda-1)a_{16}a_{34},\\
&0=N_{57}^4=-\lambda a_{18}a_{27}+a_{17} a_{28}-(1+\lambda)
a_{16}a_{38}, &
\end{align*}
we obtain $a_{16}=0$. On the other hand,
\begin{align*}
&0=N_{67}^4=(1-\lambda) a_{27} a_{28},&\quad& 0=N_{47}^2=-(1+\lambda)a_{27}a_{24},\\
&0=N_{78}^7= (1+\lambda)a_{23}a_{38}+(\lambda-1)
a_{27}a_{43},&\quad &0=N_{67}^7=2 a_{27} a_{23},\\
&0=N_{37}^8=(1-\lambda)a_{23}a_{34}-(1+\lambda)a_{27}a_{83},
\end{align*}
and \eqref{r3l44} imply that $a_{27}=0$; and from the equations
\begin{align*}
&0=N_{78}^7=(1+\lambda)a_{23}a_{38},\qquad 0=N_{68}^7=(1+\lambda)a_{23}a_{28},\\
&0=N_{47}^7=(\lambda-1)a_{23}a_{34},\qquad 0=N_{46}^7=(\lambda-1)a_{23}a_{24}
\end{align*}
we conclude that $a_{23}=0$ using again \eqref{r3l44}. Moreover,
\begin{align*}
&0=N_{48}^5=\lambda a_{21}a_{24},\qquad 0=N_{47}^5=a_{24}a_{31}+(\lambda-1)a_{21}
a_{34},\\
&0=N_{68}^5=\lambda a_{21}a_{28},\qquad
0=N_{34}^5=a_{24}a_{71}+(1+\lambda)a_{21}a_{83},
\end{align*}
imply that $a_{21}=0$ using again \eqref{r3l44}.
So, according to 
\eqref{diag},  $a_{24}{}^2+a_{28}{}^2\not=0$. Since
$\lambda\not=0,\pm1$, the equations
\begin{align*}
&0=N_{37}^6=2\lambda a_{24} a_{28},\\
&0=N_{34}^4 = (\lambda-1) a_{24} a_{43}+(1+\lambda) a_{28}a_{83},\\
&0=N_{78}^8=(\lambda-1) a_{28} a_{34}+(1+\lambda) a_{24}a_{38},
\end{align*}
imply that $a_{34}a_{43}=a_{38}a_{83}=0$. Therefore,
$$
0=N_{37}^6= 1+a_{33}{}^2+\lambda(a_{34}a_{43}-a_{38}a_{83}),
$$
which implies that $0=1+a_{33}{}^2$. But  this is not possible, and hence
$\mathcal J$ cannot be integrable on $T^*(\mathbb R\times\mathfrak r_{3,\lambda})$.
\end{proof}


\begin{prop} \label{Lie3}
The Lie algebra $\mathfrak r_{4}$ has no generalized complex structure.
\end{prop}

\begin{proof}
From  \eqref{base},  \eqref{cincofam} and \eqref{courant}, we have
that with respect to the basis $\{X_i\}$ the only non-zero  Lie
brackets on $T^*(\mathfrak r_{4})$ are
\[\begin{array}{lll}
[X_1,X_2]=X_2,&[X_1,X_3]=X_2+X_3,&[X_1,X_4]=X_3+X_4,\\[5pt]
[X_1,X_7]=-X_7- X_8,&[X_1,X_6]=-X_6-X_7, &[X_1,X_8]=-X_8, \\[5pt]
[X_2,X_6]=[X_3,X_6]=X_5,&[X_3,X_7]=[X_4,X_7]=X_5,&[X_4,X_8]=X_5.
\end{array}
\]
If there is a hermitian structure
$(\mathcal J,\langle\cdot\, ,\cdot\rangle)$
on
$T^*(\mathfrak r_{4})$, then
\begin{align*}
&0=N_{78}^1=a_{18}{}^2,\qquad 0=N_{67}^1=a_{17}{}^2-a_{16}a_{18},\\
&0=N_{35}^6=a_{12}{}^2,\qquad 0=N_{35}^8=-a_{13}{}^2+a_{12}a_{14},
\end{align*}
imply that
$a_{12}=a_{13}=a_{17}=a_{18}=0$. So,
it must be $a_{16}{}^2+a_{14}{}^2\not=0$
since, according to \eqref{diag},  at least an element $a_{1j}$ of the first row
of the matrix \eqref{jota88} associated to $\mathcal J$ must be non-zero
for $j\not=1$. Now, we
consider
\begin{align*}
&0=N_{78}^2=2 a_{16} a_{38},& &0=N_{78}^8= 2 a_{14} a_{38}\\
&0=N_{26}^7= -2a_{16}a_{72},&  &0=N_{24}^7=-2 a_{14}a_{72}\\
&0=N_{47}^7=a_{14}(a_{43}-a_{32}),&   &0=N_{67}^7= a_{16}(a_{43}-a_{32})\\
&0=N_{15}^3= a_{14}a_{38}-a_{16}a_{32},& &0=N_{35}^5=
a_{14}a_{43}+a_{16} a_{72}.
\end{align*}
From these equations, and using that $a_{16}$ and $a_{14}$ do
not vanish simultaneously, we
conclude that $a_{38}=a_{72}=a_{43}=a_{32}=0$. Then,
\[
0=N_{37}^5 = 1+a_{33}{}^2+ a_{32}(a_{11}+a_{33}+a_{23})+a_{43}(a_{33}-a_{11}+a_{34})-
a_{27}a_{72}-a_{83}a_{38},
\]
implies that $0=1+a_{33}{}^2$, which is not possible. Thus,
$\mathcal J$ cannot be integrable.
\end{proof}

\begin{prop} \label{Lie4}
For $\lambda\not=0,\pm1$, the Lie algebra $\mathfrak r_{4,\lambda}$ does not admit
generalized complex structures.
\end{prop}
\begin{proof}
With respect to the basis $\{X_i\}$ given by  ~\eqref{base}, and
according to~\eqref{courant} and
\eqref{cincofam}, the only non-zero Lie brackets on $T^*(\mathfrak r_{4,\lambda})$
are
\[
\begin{array}{lll}
[X_1,X_2]=X_2,&[X_1,X_3]=\lambda X_3,&[X_1,X_4]=X_3+\lambda X_4,\\[5pt]
[X_1,X_6]=-X_6,&  [X_1,X_7]=-\lambda X_7- X_8, &[X_1,X_8]=-\lambda X_8,\\[5pt]
[X_2,X_6]=[X_4,X_7]=X_5,&  [X_3,X_7]=[X_4,X_8]=\lambda X_5.
\end{array}
\]
Suppose that, for
$\lambda\in\mathbb R-\{-1,0,1\}$,  $T^*(\mathfrak r_{4,\lambda})$ has a
hermitian structure with complex structure $\mathcal J$.
Since all the coefficients $N^k_{ij}$ of the Nijenhuis tensor
of $\mathcal J$ are zero, we have
\[
0=N_{78}^1=a_{18}{}^2,\qquad \mbox{and}\quad 0=N_{45}^7= a_{13}{}^2 .
\]
Thus, $a_{18}=a_{13}=0$. Let us consider the equation
\begin{equation}\label{cuad22}
\begin{split}
0=N_{26}^6=&\ \ 1+a_{22}{}^2+a_{12}a_{12}+a_{16}a_{61}+a_{23}a_{42}-a_{28}a_{72}\\&
\ \ +\lambda(a_{23}a_{32}+a_{24}a_{42}-a_{27}a_{72}-a_{28}a_{82}).
\end{split}
\end{equation}
Since
\begin{align*}
&0=N_{57}^2 = (\lambda-1)a_{16} a_{17},&  &0=N_{67}^4= 2\lambda
a_{16}a_{38}-(1+\lambda)a_{17}a_{28},\\
&0=N_{67}^7= a_{16}a_{43}-(1-\lambda)a_{17}a_{23} ,& &0=N_{57}^6=
(1+\lambda)a_{12}a_{17}, \\
&0=N_{37}^6=(1+\lambda)a_{17}a_{72}-a_{12}a_{43},&
 &0=N_{78}^6=2\lambda a_{12}a_{38}-(1-\lambda)a_{17}a_{42},
\end{align*}
we obtain that $a_{17}=0$ using \eqref{cuad22} and the conditions $\lambda\not=0,1,-1$. Now,
the equations
\begin{align*}
&0=N_{46}^1=-(1+\lambda)a_{14}a_{16},&
&0=N_{46}^4=(1+\lambda)a_{14}a_{28}-a_{16}a_{43},\\
&0=N_{46}^7=(1-\lambda)a_{14}a_{23}+2\lambda a_{16}a_{83},&
&0=N_{45}^6=(\lambda-1)a_{14}a_{12},\\
&0=N_{34}^6=(1+\lambda)a_{14}a_{72}+2\lambda a_{12}a_{83},&
&0=N_{25}^5=a_{12}(a_{11}+a_{22})+\lambda a_{14}a_{42}
\end{align*}
and \eqref{cuad22} imply that $a_{14}=0$. Hence $a_{12}{}^2+a_{16}{}^2\not=0$.
Morever, we have
\begin{align*}
&0=N_{37}^2= -a_{16}a_{43},& &0=N_{37}^6=-a_{12}a_{43},&
&0=N_{78}^2=2\lambda a_{16}a_{38},\\
&0=N_{78}^6=2\lambda a_{12}a_{38},& &0=N_{56}^6=2a_{12}a_{16},&
&0=N_{35}^5=a_{12}a_{23}+a_{16}a_{72},\\
&0=N_{15}^3=-a_{12}a_{27}-a_{16}a_{32}.
\end{align*}
Therefore, $a_{43}=a_{38}=0$ and
$a_{23}a_{32}=a_{27}a_{72}=0$. Taking into account these
equalities and
$$
0=N_{37}^5=\lambda(1+a_{33}{}^2)+a_{23}a_{32}-a_{27}a_{72}+a_{43}(a_{33}-a_{11}+\lambda
a_{34})-\lambda a_{38}a_{83}
$$
we conclude that $1+a_{33}{}^2=0$ because $\lambda \not=0$. This
proves that, for $\lambda\not=-1,0,1$, $\mathfrak r_{4,\lambda}$
does not admit  generalized complex structures.
\end{proof}

\begin{prop}\label{Lie5}
The Lie algebra $\mathfrak r_{4,\mu,\lambda}$ has no generalized complex structure
for   \\ $-1<\mu<\lambda<1$ such that $\mu\lambda\not=0$ and $\mu+\lambda\not=0$.
\end{prop}
\begin{proof}
For $T^*(\mathfrak r_{4,\mu,\lambda})$ we take the basis $\{X_i\}$
defined by ~\eqref{base}. Then, using~\eqref{cincofam} and \eqref{courant} we see
that the only non-zero Lie brackets on $T^*(\mathfrak r_{4,\mu,\lambda})$ are
\[
\begin{array}{lll}
[X_1,X_2]=X_2,&[X_1,X_3]=\mu X_3,&[X_1,X_4]=\lambda X_4,\\[5pt]
[X_2,X_6]=X_5,&[X_3,X_7]=\mu X_5,&[X_4,X_8]=\lambda X_5,\\[5pt]
[X_1,X_6]=-X_6,&[X_1,X_7]=-\mu X_7,&[X_1,X_8]=-\lambda X_8.
\end{array}
\]
As in proof of the previous propositions, we assume that
$T^*(\mathfrak r_{4,\mu,\lambda})$ has a hermitian structure
$(\mathcal J,\langle\cdot\, ,\cdot\rangle)$. Then,
\begin{equation}\label{r4lm22}
0=N_{26}^5=
1+a_{22}{}^2+a_{21}a_{12}+a_{16}a_{61}+\mu(a_{23}a_{32}-a_{27}a_{72})+\lambda
(a_{24}a_{42}-a_{28}a_{82}).
\end{equation}
Consider
\begin{align*}
&0=N_{68}^1= (1-\lambda)a_{16}a_{18},&  &0=N_{58}^6=(1+\lambda)a_{12}a_{18},\\
&0=N_{58}^3= (\lambda-\mu)a_{17}a_{18},& &0=N_{58}^8= 2\lambda a_{14}a_{18}, \\
&0=N_{48}^6= (\lambda-1) a_{14}a_{42}+(1+\lambda)a_{18}a_{82},&  &
0=N_{48}^2=(1-\lambda)a_{18}a_{24}-(1+\lambda)a_{14}a_{28},
\end{align*}
\begin{eqnarray*}
 0&=&N_{78}^2= (1+\mu)
a_{18}a_{27}+(\lambda+\mu)a_{16}a_{38}-(1+\lambda)a_{17}a_{28},\\
0&=&N_{78}^6=(1-\mu)a_{18}a_{32}+(\lambda+\mu)a_{12}a_{38}+(\lambda-1)a_{17}a_{42}.
\end{eqnarray*}
From these equations and \eqref{r4lm22}, we obtain $a_{18}=0$. Now,
we have
\begin{align*}
&0=N_{67}^1= (1-\mu)a_{16}a_{17},& & 0=N_{57}^6= (1+\mu)a_{12}a_{17},\\
&0=N_{68}^3= (1+\lambda) a_{17}a_{28}-(\lambda+\mu)a_{16}a_{38}, &
&0=N_{78}^6= (\lambda+\mu) a_{12}a_{38}+(\lambda-1)a_{17}a_{42},\\
&0=N_{37}^1= -2\mu a_{17}a_{13},&
&0=N_{37}^6=(\mu-1)a_{13}a_{32}+(1+\mu)a_{17}a_{72}, \\
&0=N_{37}^2=(1-\mu)a_{17}a_{23}-(1+\mu)a_{13}a_{27}. &
\end{align*}
Then, using \eqref{r4lm22}, we see that $a_{17}=0$. From the equations
\begin{align*}
&0=N_{68}^8= (1+\lambda)a_{14}a_{28},&
&0=N_{48}^6=(\lambda-1)a_{14} a_{42},\\
&0=N_{56}^8= (1+\lambda)a_{14}a_{16},&
&0=N_{47}^2= -(1+\mu)a_{14}a_{27}+(\lambda-\mu)a_{16}a_{34},\\
& 0=N_{25}^8=(\lambda-1) a_{12} a_{14},&
& 0=N_{27}^8=(1-\mu)a_{14}a_{32}+(\mu-\lambda)a_{12}a_{34},
\end{align*}
and \eqref{r4lm22}, we conclude that $a_{14}=0$. Moreover, we have
\begin{align*}
&0=N_{37}^2= -(1+\mu)a_{13}a_{27},& & 0=N_{37}^6=(\mu-1)a_{13}a_{32},\\
&0=N_{36}^1=-(1+\mu)a_{13}a_{16},& &
0=N_{36}^4=(1+\lambda)a_{13}a_{28}+(\lambda-
\mu)a_{16}a_{43},\\
&0=N_{35}^6=(\mu-1)a_{12}a_{13},& &
0=N_{36}^8=(1-\lambda)a_{13}a_{24}-(\lambda+\mu) a_{16}a_{83},
\end{align*}
which  imply that $a_{13}=0$ using again \eqref{r4lm22}. Thus,
$a_{12}{}^2+a_{16}{}^2\not=0$. Now, taking account  the equations
\begin{align*}
&0=N_{56}^6=2a_{12}a_{16},& & \\
&0=N_{78}^2 = (\lambda+\mu) a_{16}a_{38},& & 0=N_{78}^6= (\lambda+\mu)a_{12}a_{38},\\
&0=N_{38}^2=(\mu-\lambda)a_{16}a_{43},& & 0=N_{38}^6=(\mu-\lambda)a_{12}a_{43},\\
&0=N_{15}^6=-a_{12}a_{24}-a_{16}a_{82},& &0=N_{15}^4=-a_{12}a_{28}-a_{16}a_{42},
\end{align*}
we have that $a_{38}=a_{43}=0$ and $a_{24}a_{42}=a_{28}a_{82}=0$. So,
$$
0=N_{48}^5=\lambda(1+a_{44}{}^2)+\mu(a_{34}a_{43}-a_{38}a_{83})+a_{24}a_{42}-a_{28}a_{82}
=\lambda(1+a_{44}{}^2).
$$
This implies that $\lambda=0$ or $1+a_{44}{}^2=0$, which is not possible.
 This completes the proof.
\end{proof}

 Let  $\frak g$ be an arbitrary Lie algebra. Denote by
$b_i(\frak g)$ the dimension of the $i$-th
cohomology group $H^{i}(\frak g)$
of $\frak g$, by $\frak g ' =[\frak g , \frak g]$
the derived subalgebra and by $\frak z (\frak g)$
the center of $\frak g$. We recall that $\frak g$ is called
\emph{completely solvable} when $\frak g$ is solvable and
ad$(x)$ has only real eigenvalues for any $x\in \frak g$.

\begin{thm}\label{main}
Let $G$ be  a four dimensional solvable  Lie group with Lie
algebra $\frak g$. Then the following statements are equivalent:
\begin{enumerate}
\item[(i)] $G$ has no left invariant generalized complex structure;
\item[(ii)]  $G$ admits neither  left invariant symplectic nor  left
invariant complex structures;
\item[(iii)] $\frak g$ is completely solvable and one of the two following
conditions is satisfied:
\begin{enumerate}
\item[(a)] $b_1(T^*\frak g)=3, \;\; b_3(T^*\frak g)=5,
\;\; b_1(\frak g /\frak z (\frak g' ))=1$, or
\item[(b)] $b_1(T^*\frak g)=1, \;\; b_3(T^*\frak g)=2$.
\end{enumerate}
\end{enumerate}
\end{thm}

\begin{proof}
Clearly (i) implies (ii). The converse follows
from Propositions~\ref{Lie1} to \ref{Lie5}. 
The calculation of the numbers $b_i(T^*(\mathfrak g))$, $(i=1,3)$,
and $b_1(\mathfrak g /\mathfrak z (\mathfrak g' ))$, where $\mathfrak g$ is a  four
dimensional solvable Lie algebra, shows that condition (iii)
is satisfied if and only if~(i) holds (see Table in Section ~\ref{sec6}).
\end{proof}


\section{Generalized complex structures of type $1$ on solvable Lie groups}\label{sec-5}


In this section we exhibit the four dimensional solvable Lie
algebras which have a generalized complex structure of type $1$;
Theorem~\ref{main} implies that they admit either symplectic or
complex structures. We also study necessary and sufficient
conditions on
  a four dimensional solvable Lie algebra $\frak g$ to admit generalized
  complex structures of type $1$. As a consequence of our results,
  we obtain in Corollary~\ref{nouno} a condition, involving the odd
numbers $b_i(T^*(\mathfrak g))$, for the non-existence of such
structures.

\subsection{Existence}\label{subsec-5}

First, we list below the family of Lie algebras having either
symplectic or complex structures. Such Lie algebras together with
those shown in \eqref{cincofam}, exhaust the class of four
dimensional solvable Lie algebras (see \cite{ABDO,D,Mu,Z}).
\begin{equation}\label{tipo1-si}
\begin{tabular}{rl}
$\mathfrak{aff}(\mathbb R)\!\times\!\mathfrak{aff}(\mathbb R)$ :\ & $[e_0,e_1]=e_1, \
[e_2,e_3]= e_3$;\\
$\mathfrak{aff}(\mathbb C)$ :\  &$[e_0,e_2]=e_2,\  [e_0,e_3]=e_3,\ [e_1, e_2]= e_3, \
[e_1,e_3] =- e_2$;\\
$\mathbb R\times\mathfrak e(2)$ :\  &$[e_1,e_2]= - e_3,\  [e_1,e_3]=e_2$;\\
$\mathbb R\times\mathfrak h_3$ :\  &$[e_1,e_2]=e_3$;\\
$\mathbb R\times\mathfrak r_{3, \lambda}$ :\  &$ [e_1, e_2]=e_2,\ [e_1,e_3]=\lambda
e_3,\quad\lambda \in \{-1,0,1\}$ ;\\
$\mathfrak r_{4,\lambda}$ :\  &$[e_0, e_1]=e_1, \ [e_0,e_2]=\lambda e_2,\ [e_0, e_3]=
e_2+\lambda e_3, \ \lambda \in \{-1,  0,1\}$;\\
$\mathfrak r_{4,\mu,1}$ :\  &$[e_0,e_1]=e_1,\  [e_0, e_2]= \mu e_2,\  [e_0,e_3]=
e_3, \quad  -1<\mu\leq 1,\  \mu \neq 0$;\\
$\mathfrak r_{4,\mu,\mu}$ :\  &$[e_0,e_1]=e_1,\  [e_0, e_2]=\mu e_2,\  [e_0,e_3]= \mu
e_3, \  -1<\mu<1,\  \mu\neq 0$;\\
$\mathfrak r_{4,\mu,-\mu}$ :\  &$[e_0,e_1]=e_1,\  [e_0, e_2]=\mu e_2,\  [e_0,e_3]=
-\mu e_3, \quad -1< \mu <0 $;\\
$\mathfrak r_{4,-1,\lambda}$ :\  &$[e_0,e_1]=e_1,\  [e_0, e_2]=- e_2,\  [e_0,e_3]=
\lambda e_3, \quad -1 <\lambda <0$; \\
$\mathfrak r_{4,-1,-1}$ :\  &$[e_0,e_1]=e_1,\  [e_0, e_2]=- e_2,\  [e_0,e_3]=-e_3$;
\\
$\mathbb R \times \mathfrak r'_{3,\lambda}$ :\  &$[e_1,e_2]=\lambda e_2 -e_3,\
[e_1,e_3] =e_2 +\lambda e_3, \quad \lambda > 0$;\\
$\mathfrak n_4$ :\ & $[e_0, e_1]= e_{2},\  [e_0, e_2]= e_3$;\\
$\mathfrak r'_{4,\mu,\lambda}$ :\  &$[e_0,e_1]=\mu e_1,\ [e_0,e_2]=\lambda e_2 -e_3,
\  [e_0,e_3]=e_2 +\lambda e_3,$\\
&\ \ $\mu > 0,\  \lambda \in \mathbb R$;\\
$\mathfrak d_4$ :\  &$ [e_0,e_1]=  e_1,\ [e_0,e_2]= -e_2,\ [e_1,e_2]=e_3$;\\
$\mathfrak d_{4,\lambda}$ :\  &$[e_0,e_1]= \lambda e_1,\ [e_0,e_2]= (1-\lambda)e_2,
\ [e_0,e_3]=  e_3, \  [e_1,e_2]=e_3,$ \\
& $\lambda \geq 1/2$;\\
$\mathfrak d'_{4,\lambda}$ :\  &$[e_0, e_1]=\lambda e_1-e_2,\ [e_0, e_2]=
e_1+\lambda e_2,\  [e_0, e_3]= 2\lambda e_3,$ \\ &$[e_1,e_2]=e_3, \qquad
\lambda \geq 0 $;\\
$\mathfrak h_4$ :\  &$[e_0,e_1]=e_1, \  [e_0,e_2]=e_1+e_2, \
[e_0,e_3]=2e_3, \ [e_1,e_2]=e_3$.
\end{tabular}
\end{equation}

\begin{prop} \label{todo-tipo}
The Lie algebras
$\mathfrak{aff}(\mathbb R)\times\mathfrak{aff}(\mathbb R)$,
$\mathfrak{aff}(\mathbb C)$,
$ \mathbb R \times \mathfrak e(2)$,
$\mathbb R\times\mathfrak h_3$,
$\mathbb R\times\mathfrak r_{3,0}$,
$\mathfrak r_{4,-1,-1}$,
$\mathfrak r_{4,\mu ,0}' \; (\mu >0)$,
$\mathfrak d_{4,\frac 12}$,
$\mathfrak d_{4,2}$ and
$\mathfrak d'_{4,\lambda} \; (\lambda > 0)$
admit generalized complex structures of type $0$, $1$ and $2$.
\end{prop}
\begin{proof}
It follows from results in \cite{MR,O,SJ} that all of the above Lie algebras admit
both symplectic and complex  structures, which give rise to generalized complex
structures of type $0$ and $2$, respectively.  A generalized complex structure of
type $1$ on $\;\mathfrak{aff}(\mathbb R) \times
\mathfrak{aff}(\mathbb R)\;$ can be obtained by combining one of type $0$ with one of
type $1$ on $\;\mathfrak{aff}(\mathbb R)\;$ (see \eqref{aff(R)-0} and
\eqref{aff(R)-1}).
For the remaining Lie algebras, we exhibit a generalized complex structure of
type~$1$.
\begin{equation} \begin{split}
\mathfrak{aff}(\mathbb C),\ \mathfrak r_{4,-1,-1}:\hspace{1cm} & \mathcal J(e_0)
=\alpha^1,\quad
\mathcal J(e_2)=e_3
\\
\mathbb R\times \mathfrak e(2),\ \mathbb R\times \mathfrak r_{3,0},\
\mathfrak d_{4,\frac 12},\ \mathfrak d'_{4,\lambda}\ (\lambda>0):\hspace{1cm} &
\mathcal J(e_0)=\alpha^3,\quad \mathcal J(e_1) = e_2,
\\
\mathbb R\times \mathfrak h_3,\  \mathfrak r'_{4,\mu,0}\ (\mu>0),\  \mathfrak d_{4,2}
:\hspace{1cm} &\mathcal J(e_0) = e_1,\quad \mathcal J(e_2)=\alpha^3.
\end{split}
\end{equation}
\end{proof}

\begin{prop} The Lie algebras
$\mathbb R\times\mathfrak r_{3,-1}$, $\mathfrak r_{4,-1}$,
$\mathfrak r_{4,0}$, $\mathfrak n_4$, $\mathfrak r_{4,\mu,-\mu}\;
(-1<\mu<0)$ and $\mathfrak r_{4,-1,\lambda}\;(-1\leq\lambda<0)$
admit generalized complex structures of type $0$ and $1$, but not of
type~$2$.
\end{prop}
\begin{proof} First we notice that every Lie algebra mentioned in the proposition
has symplectic structures but does
not admit complex structures (\cite{MR,O1,O,SJ}), so it does not possess generalized
complex structures of type~2 (Theorem~\ref{simple}). For each one of these Lie algebras,
 we show a generalized complex
structure of type $1$:
\begin{equation} \begin{split}
\mathbb R\times\mathfrak r_{3,-1},\ \mathfrak r_{4,0},\ \mathfrak
r_{4,\mu,-\mu},\ \mathfrak n_4:\hspace{1cm} &
\mathcal J(e_0) =e_1,\quad\mathcal J(e_2)=\alpha^3,\\
\mathfrak r_{4,-1},\ \mathfrak r_{4,-1,\lambda}:\hspace{1cm} &\mathcal J(e_0)
=e_3,\quad\mathcal J(e_1)=\alpha^2.
\end{split}
\end{equation}
\end{proof}

\begin{prop} The Lie algebras
$\mathbb R\times\mathfrak r_{3,1}$,
$\mathfrak r_{4,1}$,
$\mathfrak r'_{4,\mu,\lambda}\; (\mu>0,\;\lambda\not=0)$,
$\mathbb R\times\mathfrak  r_{3,\lambda}'\; (\lambda \neq 0)$,
$\mathfrak  r_{4,\mu,\mu}\;  (-1<\mu\leq 1, \, \mu \neq 0)$,
$\mathfrak r_{4,\mu,1}\; (-1<\mu\leq 1, \, \mu \neq 0)$, and
$\mathfrak d'_{4,0}$  admit generalized complex structures of type  $1$ and $2$, but
not of type $0$.
\end{prop}
\begin{proof} These Lie algebras have complex structures and do not
admit symplectic structures (\cite{MR,O1,O,SJ}), thus they admit generalized complex
structures of type~2 but not of type~0.
 A generalized complex structure of type~1 is given by
\begin{equation} \begin{split}
\mathbb R\times\mathfrak r_{3,1},\ \mathfrak r_{4,\mu,\mu},\ \mathfrak
r'_{4,\mu,\lambda}:\hspace{1cm}&\mathcal J(e_0)=\alpha^1,\quad \mathcal J(e_2) =
e_3,\\
\mathfrak r_{4,1},\ \mathfrak r_{4,\lambda,1}:\hspace{1cm} &\mathcal
J(e_0)=\alpha^2,\quad \mathcal J(e_1) = e_3,\\
\mathbb R\times \mathfrak r'_{3,\lambda},\ \mathfrak d'_{4,0}:\hspace{1cm}
&\mathcal J(e_0)=\alpha^3,\quad \mathcal J(e_1) = e_2,
\end{split}
\end{equation}
\end{proof}

\subsection{Obstructions }\label{sec6}


\begin{prop} \label{prop-6.1}
$\mathfrak d_4 $ is the unique four dimensional solvable Lie algebra admitting
generalized complex structures of type  $2$, but not of type $0$ or $1$.
\end{prop}
\begin{proof}
We consider the basis $\{X_i\}$ defined by \eqref{base} for $T^*(\mathfrak d_4)$.
Taking into account \eqref{courant} and the structure equations of the Lie algebra
$\mathfrak d_4$ given in \eqref{tipo1-si}, we see that the only non-zero brackets on
$T^*(\mathfrak d_4)$ are
\[\begin{array}{lll}
[X_1,X_2]= X_2,&[X_1,X_3]=-X_3,&[X_2,X_3]=X_4,\\[5pt]
[X_2,X_6]=-[X_3,X_7]=X_5,&-[X_1,X_6]=[X_3,X_7]=X_6,&[X_1,X_7]=-[X_2,X_8]=X_7.
\end{array}
\]
Supposse that $\mathcal J$ is a generalized complex structure on $\mathfrak
d_4$. Let us consider the equations
\begin{align*}
&0=N_{56}^7= a_{16} a_{24},\quad 0=N_{56}^8= a_{16}a_{14},\quad 0=N_{46}^3=-a_{34}
a_{16}-a_{24} a_{17},\\
&0=N_{26}^8=a_{24}(a_{12}-a_{34})-a_{16}a_{82},\\
&0=N_{36}^8=-2a_{14}a_{23}-a_{13}a_{24}+a_{24}{}^2+a_{16}a_{83}.
\end{align*} The condition $\mathcal J^2=-\id$ implies that
\begin{align*}
0=&\ (\mathcal
J^2)_2^4=a_{14}a_{21}+a_{24}(a_{22}+a_{44})+a_{23}a_{34}+a_{16}a_{81}-a_{27}a_{83},\\
-1=&\ (\mathcal
J^2)_4^4=a_{14}a_{41}+a_{24}a_{42}+a_{34}a_{43}+a_{44}{}^2+a_{18}a_{81}+a_{28}a_{82}
+a_{38}a_{83},
\end{align*}
and so we obtain $a_{16}=0$. Now, from the equations
\begin{align*}
&0=N_{47}^1= a_{14} a_{17},\quad 0=N_{47}^2= a_{17}a_{24},\quad 0=N_{57}^6=-a_{34}
a_{17},\\
&0=N_{47}^7=a_{34}(a_{13}-a_{24})-a_{17}a_{83},\\
&0=N_{47}^6=-2a_{14}a_{32}-a_{12}a_{34}+a_{34}{}^2+a_{17}a_{82},
\end{align*} and, from $\mathcal J^2=-\id$,
\begin{align*}
0=&\ (\mathcal J^2)_2^4=a_{14}a_{31}+a_{24}a_{32}+a_{34}(a_{33}+a_{44})+a_{17}a_{81}
+a_{27}a_{82},\\
-1=&\
(\mathcal J^2)_4^4=a_{14}a_{41}+a_{24}a_{42}+a_{34}a_{43}+a_{44}{}^2+a_{18}a_{81}
+a_{28}a_{82} +a_{38}a_{83},
\end{align*}
we obtain $a_{17}=0$. Moreover,  $a_{27}=0$ because $N_{78}^3= -a_{27}{}^2=0$.

The equations
\begin{align*}
&0=(\mathcal J^2)_1^5=-2 a_{14}a_{18}, \quad 0=N_{78}^8=a_{18}a_{34}-a_{14}a_{38},
\quad 0=N_{58}^7=a_{14}a_{28}-a_{13}a_{18},\\
&0=N_{68}^8= a_{14} a_{28}-a_{18}a_{24},\quad
0=N_{68}^7=a_{28}(a_{13}+a_{24})-a_{18}a_{23},\\
&0=N_{17}^7=-1+a_{13}a_{31}-a_{23}a_{32}+a_{33}{}^2-a_{21}a_{34},
\end{align*}
imply that $a_{18}=0$. But, since
\begin{align*}
& 0=(\mathcal J^2)_2^6=-2a_{24}a_{28},
\quad 0=N_{68}^7=a_{28}(a_{13}+a_{24}),\\
&0=(\mathcal J^2)_1^6=-a_{14}a_{28}, \quad
0=(\mathcal J^2)_1^8=a_{12}a_{28}+a_{13}a_{38}, \\
&-1=(\mathcal J^2)_1^1=a_{11}{}^2+a_{12}a_{21}+a_{13}a_{31}+a_{14}a_{41},
\end{align*}
we obtain $a_{28}=0$. And finally, from the equations
\begin{align*}
& 0=(\mathcal J^2)_1^8=a_{23}a_{38},
\quad 0=N_{18}^2=-a_{38}a_{13},\quad 0=(\mathcal J^2)_3^7=-2a_{34}a_{38}, \\
&0=N_{17}^7=1+a_{33}{}^2+a_{13}a_{31}-a_{23}a_{32}-a_{21}a_{34},
\end{align*}
we have $a_{38}=0$. So the matrix $\mathcal J_2$ in \eqref{jotas4} is the null
matrix, and then $\mathfrak d_4$ does not admit generalized complex structures of
types 0 and 1. The almost complex structure defined by $J(e_0) = e_1$ and $J(e_3) =
e_2$ is integrable and thus $\mathfrak d_4$ admits a generalized complex structure
of type 2. The uniqueness is seen in the table at the end of this section.
\end{proof}

\begin{prop} \label{prop-6.2} The Lie algebras $\mathfrak d_{4,\lambda}$, $(\lambda \neq \frac
12 ,\,2)$ and $\mathfrak h_4$ admit generalized complex structures of type
$0$ and $2$, but not of type $1$.
\end{prop}
\begin{proof} Doing a  similar calculation to
that made in the previous proposition, one can check that
 for a generalized complex structure on
$\mathfrak h_4$, the matrices
$\mathcal J_1$ and $\mathcal J_2$ in \eqref{jotas4}
are
\[
\mathcal J_1=\begin{pmatrix} a_{11}&0&a_{13}&0\\ a_{21}&a_{22}&a_{23}&a_{13}\\
a_{31}&0&a_{33}&0\\ a_{41}&a_{42}&a_{43}&a_{44}\end{pmatrix},
\qquad
\mathcal J_2=\begin{pmatrix} 0&0&0&a_{18}\\ 0&0&2a_{18}&a_{28}\\
0&-2a_{18}&0&a_{38}\\ -a_{18}&-a_{28}&-a_{38}&0\end{pmatrix}.
\]
From $\mathcal J^2=-\id$ we have $a_{18}{}^2+a_{13}^2\not=0$.
Since $0=N_{68}^8= 3 a_{18}a_{13}$,  we have the two following
possibilities:
\begin{itemize}
\item\ \ $a_{18}\not=0$, $a_{13}=0$. Then, $\rank\mathcal J_2=4$
and the possible generalized complex structures are of type 0; for example,
\[
\mathcal J(e_0)=2 \alpha^3,\qquad \mathcal J(e_1)= \alpha^2.
\]
\item\ \  $a_{13}\not=0$, $a_{18}=0$. In this case, we obtain
$a_{38}=a_{28}=0$ using
$0=N_{78}^7= 3 a_{13}a_{38}$ and $0=N_{68}^7=a_{13}(4 a_{28}+a_{38})$.
So $\mathcal J_2\equiv 0$ and the possible
generalized complex structures are of type 2; for example,
\[
\mathcal J(e_0)=e_2,\qquad \mathcal J(e_1)= e_3.
\]
\end{itemize}

For a generalized complex structure on
the Lie algebra $\mathfrak  d_{4, \lambda}$, $(\lambda \neq \frac
12 , \,1,\, 2)$, the matrices $\mathcal J_1$ and $\mathcal J_2$, given by
\eqref{jotas4}, are
\[
\mathcal J_1=\begin{pmatrix} a_{11}&a_{12}&a_{1,3}&0\\ a_{21}&a_{22}&a_{23}&a_{24}\\
a_{31}&a_{32}&a_{33}&a_{34}\\ a_{41}&a_{42}&a_{43}&a_{44}\end{pmatrix},\qquad
\mathcal J_2=\begin{pmatrix} 0&0&0&a_{18}\\ 0&0&a_{18}&a_{28}\\
0&-a_{18}&0&a_{38}\\ -a_{18}&-a_{28}&-a_{38}&0\end{pmatrix}.
\]
We consider two possibilities according to $a_{18}\not=0$ or $a_{18}=0$:
\begin{itemize}
\item\ \ If $a_{18}\not=0$, $\rank\mathcal J_2=4$ and the
generalized complex structures are of type 0; for example,
\[
\mathcal J(e_0) = \alpha^3,\qquad \mathcal J(e_1) =\alpha^2.
\]
\item\ \ If $a_{18}=0$, then $a_{12}{}^2+a_{13}{}^2\not=0$. Since $0=
N_{35}^6= (1-2\lambda)a_{12}a_{13}$, we consider two subcases
\begin{itemize}
\item[A)]\ \ $a_{12}\not=0$, $a_{13}=0$. From $N_{58}^5=-\lambda a_{12}a_{28}$ we
obtain $a_{28}=0$ and from $0=N_{78}^6=-a_{38}(a_{34}-a_{12}(\lambda-2))$ and $0=
N_{47}^6= a_{34}(\lambda a_{12}+a_{34})$ we obtain $a_{38}=0$. Hence $\mathcal
J_2\equiv 0$ and the generalized complex structures are of type 2; for example,
\[
\mathcal J(e_0) = \lambda e_1,\qquad \mathcal J(e_2) =-e_3.
\]
\item[B)]\ \ $a_{13}\not=0$, $a_{12}=0$. From $N_{58}^5=(-1+\lambda) a_{13}a_{38}$ we
obtain $a_{38}=0$ and from $0=N_{68}^7=a_{28}(a_{24}+a_{13}(\lambda+1))$ and $0=
N_{46}^7= -a_{24}(a_{24}+a_{13}(\lambda-1))$ we obtain $a_{28}=0$. So, $\mathcal
J_2\equiv 0$ and the generalized complex structures are of type 2; for example,
\[
\mathcal J(e_0) =(1-\lambda) e_2,\qquad \mathcal J(e_1) =e_3.
\]
\end{itemize}
\end{itemize}

For a generalized complex structure on
the Lie algebra $\mathfrak  d_{4,1}$, the matrices $\mathcal J_1$ and
$\mathcal J_2$ in  \eqref{jotas4}
are
\[
\mathcal J_1=\begin{pmatrix} a_{11}&a_{12}&0&0\\ a_{21}&a_{22}&a_{23}&0\\
a_{31}&a_{32}&a_{33}&-a_{12}\\ a_{41}&a_{42}&a_{43}&a_{44}\end{pmatrix},\qquad
\mathcal J_2=\begin{pmatrix} 0&0&0&a_{18}\\ 0&0&a_{18}&a_{28}\\
0&-a_{18}&0&a_{38}\\ -a_{18}&-a_{28}&-a_{38}&0\end{pmatrix}.
\]
Therefore,  $a_{12}{}^2+a_{18}{}^2\not=0$. Since $0=N_{58}^6= 2
a_{12}a_{18}$ we consider the two following possibilities:
\begin{itemize}
\item\ \ $a_{18}\not=0$, $a_{12}=0$. Then, $\rank \mathcal J_2=4$ and
the generalized complex structures are of type 0; for example,
\[
\mathcal J(e_0) = \alpha^3,\qquad \mathcal J(e_1) = \alpha^2.
\]
\item\ \ $a_{18}=0$, $a_{12}\not=0$.  Because $N_{78}^7=-a_{12}a_{28}$ and
$0=N_{78}^6=2a_{38}a_{12}$, we have $a_{28}=a_{38}=0$. So, $\mathcal
J_2\equiv 0$ and the generalized complex structures are of type 2; for example,
\[
\mathcal J(e_0) =e_1,\qquad \mathcal J(e_2) =-e_3.
\]
\end{itemize}
\end{proof}

The previous propositions together with Table~\ref{tabla} imply
the next result.

\begin{cor}\label{nouno}
 Let $\mathfrak  g$ be a  four dimensional Lie algebra admitting a generalized complex
structure.  Then, $\mathfrak  g$ does not admit a generalized complex structure of
type $1$ if and only if $\mathfrak  g$ is completely solvable and one of the
following conditions is satisfied:
\begin{enumerate}
\item[(i)] $b_1(T^*\mathfrak  g )=b_3(T^*\mathfrak  g )=1,$ or
\item[(ii)] $b_1(T^*\mathfrak  g )=2,\;\; b_3(T^*\mathfrak  g )=4$.
\end{enumerate}
\end{cor}

\medskip
\begin{remark}
We must notice that the Lie algebra $\mathfrak d'_{4,\lambda}$
satisfies
$b_1(T^*\mathfrak d'_{4,\lambda})=
  b_3(T^*\mathfrak d'_{4,\lambda})=1$, but it
is not completely solvable. Therefore,
according to the previous Corollary,
  it has generalized complex structures of type $1$.
  In general, in the table below, the Lie algebras
  $\mathfrak g'$
are not completely solvable, so
   they have generalized complex structures of type $1$.
  \end{remark}

\medskip
In the table below we summarize the previous results and, for each
solvable Lie algebra admitting generalized complex structures, we
exhibit  \emph{one} of the simplest examples of each type (---
stands for non existence).

\setlongtables
\renewcommand{\arraystretch}{1.5}
\begin{longtable}{|p{3.7cm}|c|c|p{2.4cm}|p{2.4cm}|p{2.4cm}|}
 \hline \centerline{$\mathfrak g $}&
$b_1(T^*\mathfrak g)$ &$b_3(T^*\mathfrak g)$ & \centerline{Type
0}& \centerline{Type 1}&
\centerline{Type 2}\\
\hline
\endhead
\endfoot
\hline
\centerline{$\mathbb R\times \mathfrak r_3$}&3 &5 &
\centerline{---}&
\centerline{---}&
\centerline{---}\\
\hline
\centerline{$\mathbb R\times \mathfrak r_{3,\lambda}$}
\centerline{$|\lambda|<1,\ \lambda\not=0$}&3 &5 &
\centerline{---}&
\centerline{---}&
\centerline{---}\\
\hline
\centerline{$\mathfrak r_4$}&1&2&
\centerline{---}&
\centerline{---}&
\centerline{---}\\
\hline
\centerline{$\mathfrak r_{4,\lambda}$}
\centerline{ $\lambda\in\mathbb R,\ \lambda\not=-1,0,1$}&1&2&
\centerline{---}&
\centerline{---}&
\centerline{---}\\
\hline
\centerline{$\mathfrak r_{4,\mu,\lambda}$}
\centerline{$-1<\mu<\lambda<1$}
\centerline{$\mu\lambda\not=0,\ \mu+\lambda\not=0$}&1&2&
\centerline{---}&
\centerline{---}&
\centerline{---}\\
\hline
\centerline{$\mathfrak{aff}(\mathbb R)\times\mathfrak{aff}(\mathbb R)$}&3&5&
\centerline{$\mathcal J(e_0)=\alpha^1$,}
\centerline{$\mathcal J(e_2)=\alpha^3$}&
\centerline{$\mathcal J(e_0)=\alpha^1$,}
\centerline{$\mathcal J(e_2)=e_3$}&
\centerline{$\mathcal J(e_0)=e_1$,}
\centerline{$\mathcal J(e_2)=e_3$}\\
\hline
\centerline{$\mathfrak {aff}(\mathbb C)$}& 2 & 2 &
\centerline{$\mathcal J(e_0)=\alpha^3$,}
\centerline{$\mathcal J(e_1)=\alpha^2$}&
\centerline{$\mathcal J(e_0)=\alpha^1$,}
\centerline{$\mathcal J(e_2)=e_3$}&
\centerline{$\mathcal J(e_0)=e_3$,}
\centerline{$\mathcal J(e_1)=-e_2$}\\
\hline
\centerline{$\mathbb R\times\mathfrak h_3$}& 5 & 31 &
\centerline{$\mathcal J(e_0)=\alpha^2$,}
\centerline{$\mathcal J(e_1)=\alpha^3$}&
\centerline{$\mathcal J(e_0)=e_1$,}
\centerline{$\mathcal J(e_2)=\alpha^3$}&
\centerline{$\mathcal J(e_0)=e_1$,}
\centerline{$\mathcal J(e_2)=e_3$}\\
\hline
\centerline{$\mathbb R\times\mathfrak r_{3,-1}$}& 3 & 13 &
\centerline{$\mathcal J(e_0)=\alpha^2$,}
\centerline{$\mathcal J(e_1)=-\alpha^3$}&
\centerline{$\mathcal J(e_0)=e_1$,}
\centerline{$\mathcal J(e_2)=\alpha^3$}&
\centerline{---}\\
\hline
\centerline{$\mathbb R\times\mathfrak r_{3,0}$}& 5 & 11 &
\centerline{$\mathcal J(e_0)=\alpha^2$,}
\centerline{$\mathcal J(e_1)=-\alpha^3$}&
\centerline{$\mathcal J(e_0)=\alpha^3$,}
\centerline{$\mathcal J(e_1)=e_2$}&
\centerline{$\mathcal J(e_0)=e_3$,}
\centerline{$\mathcal J(e_1)=e_2$}\\
\hline
\centerline{$\mathbb R\times\mathfrak r_{3,1}$}& 3 & 13 &
\centerline{---}&
\centerline{$\mathcal J(e_0)=\alpha^1$,}
\centerline{$\mathcal J(e_2)=e_3$}&
\centerline{$\mathcal J(e_0)=e_1$,}
\centerline{$\mathcal J(e_2)=e_3$}\\
\hline
\centerline{$\mathfrak r_{4,-1}$}& 1 & 4 &
\centerline{$\mathcal J(e_0)=\alpha^2$,}
\centerline{$\mathcal J(e_1)=\alpha^3$}&
\centerline{$\mathcal J(e_0)=e_3$,}
\centerline{$\mathcal J(e_1)=\alpha^2$}&
\centerline{---}\\
\hline
\centerline{$\mathfrak r_{4,0}$}& 3 & 7 &
\centerline{$\mathcal J(e_0)=\alpha^1$,}
\centerline{$\mathcal J(e_2)=\alpha^3$}&
\centerline{$\mathcal J(e_0)=e_1$,}
\centerline{$\mathcal J(e_2)=\alpha^3$}&
\centerline{---}\\
\hline
\centerline{$\mathfrak r_{4,1}$}& 1 & 4 &
\centerline{---}&
\centerline{$\mathcal J(e_0)=\alpha^2$,}
\centerline{$\mathcal J(e_1)=e_3$}&
\centerline{$\mathcal J(e_0)=e_3$,}
\centerline{$\mathcal J(e_1)=e_2$}\\
\hline
\centerline{$\mathfrak r_{4,\mu,1}$,}
\centerline{$-1<\mu\leq1$, $\mu\not=0$}&1&4&
\centerline{---}&
\centerline{$\mathcal J(e_0)=\alpha^2$,}
\centerline{$\mathcal J(e_1)=e_3$}&
\centerline{$\mathcal J(e_0)=e_2$,}
\centerline{$\mathcal J(e_1)=e_3$}\\
\hline
\centerline{$\mathfrak r_{4,\mu,\mu}$,}
\centerline{$-1<\mu\leq1$, $\mu\not=0$}&1&4&
\centerline{---}&
\centerline{$\mathcal J(e_0)=\alpha^1$,}
\centerline{$\mathcal J(e_2)=e_3$}&
\centerline{$\mathcal J(e_0)=e_1$,}
\centerline{$\mathcal J(e_2)=e_3$}\\
\hline
\centerline{$\mathfrak r_{4,-1,\lambda}$,}
\centerline{$-1<\lambda<0$}&1&4&
\centerline{$\mathcal J(e_0)=\alpha^3$,}
\centerline{$\mathcal J(e_1)=\alpha^2$}&
\centerline{$\mathcal J(e_0)=e_3$,}
\centerline{ $\mathcal J(e_1)=\alpha^2$}&
\centerline{---}\\
\hline
\centerline{$\mathfrak r_{4,-1,-1}$,}&1&4&
\centerline{$\mathcal J(e_0)=\alpha^3$,}
\centerline{$\mathcal J(e_1)=\alpha^2$}&
\centerline{$\mathcal J(e_0)=\alpha^1$,}
\centerline{ $\mathcal J(e_2)=e_3$}&
\centerline{$\mathcal J(e_0)=e_1$,}
\centerline{ $\mathcal J(e_2)=e_3$}\\
\hline \centerline{$\mathbb R\times \mathfrak r'_{3,0}$}&5&11&
\centerline{$\mathcal J(e_0)=\alpha^3$,} \centerline{$\mathcal
J(e_1)=\alpha^2$}& \centerline{$\mathcal J(e_0)=\alpha^3$,}
\centerline{$\mathcal J(e_1)=e_2$}& \centerline{$\mathcal
J(e_0)=e_3$,}
\centerline{$\mathcal J(e_1)=e_2$}\\
\hline
\centerline{$\mathbb R\times \mathfrak r'_{3,\lambda}$,}
\centerline{$\lambda>0$}&3&5&
\centerline{---}&
\centerline{$\mathcal J(e_0)=\alpha^3$,}
\centerline{$\mathcal J(e_1)=e_2$}&
\centerline{$\mathcal J(e_0)=e_3$,}
\centerline{$\mathcal J(e_1)=e_2$}\\
\hline
\centerline{$\mathfrak n_4$}&3&14&
\centerline{$\mathcal J(e_0)=\alpha^3$,}
\centerline{$\mathcal J(e_1)=\alpha^2$}&
\centerline{$\mathcal J(e_0)=e_1$,}
\centerline{$\mathcal J(e_2)=\alpha^3$}&
\centerline{---}\\
\hline \centerline{$\mathfrak r'_{4,\mu,0}$,}
\centerline{$\mu>0$}&1&4& \centerline{$\mathcal J(e_0)=\alpha^1$,}
\centerline{ $\mathcal J(e_2)=\alpha^3$}& \centerline{$\mathcal
J(e_0)=e_1$,} \centerline{$\mathcal J(e_2)=\alpha ^3$}&
\centerline{$\mathcal J(e_0)=e_1$,}
\centerline{$\mathcal J(e_2)=e_3$}\\
\hline
\centerline{$\mathfrak r'_{4,\mu,\lambda}$,}
\centerline{$\mu>0$, $\lambda\not=0$}&1&2&
\centerline{---}&
\centerline{$\mathcal J(e_0)=\alpha^1$,}
\centerline{$\mathcal J(e_2)=e_3$}&
\centerline{$\mathcal J(e_0)=e_1$,}
\centerline{$\mathcal J(e_2)=e_3$}\\
\hline \centerline{$\mathfrak h_4$}&1&1& \centerline{$\mathcal
J(e_0)=2\alpha^3$,} \centerline{$\mathcal J(e_1)=\alpha^2$}&
\centerline{---}& \centerline{$\mathcal J(e_0)=e_2$,}
\centerline{$\mathcal J(e_1)=e_3$}\\
\hline
\centerline{$\mathfrak d_4$}&2&4&
\centerline{---}&
\centerline{---}&
\centerline{$\mathcal J(e_0)=e_1$,}
\centerline{$\mathcal J(e_3)=e_2$}\\
\hline
\centerline{$\mathfrak d_{4,\lambda}$}
\centerline{$\lambda>\frac 12$, $\lambda\not=1,2$}&1&1&
\centerline{$\mathcal J(e_0)=\alpha^3$,}
\centerline{$\mathcal J(e_1)=\alpha^2$}&
\centerline{---}&
\centerline{$\mathcal J(e_0)=\lambda e_1$,}
\centerline{$\mathcal J(e_2)=-e_3$}\\
\hline
\centerline{$\mathfrak d_{4,1}$}&2&4&
\centerline{$\mathcal J(e_0)=\alpha^3$,}
\centerline{$\mathcal J(e_1)=\alpha^2$}&
\centerline{---}&
\centerline{$\mathcal J(e_0)=e_1$,}
\centerline{$\mathcal J(e_2)=-e_3$}\\
\hline
\centerline{$\mathfrak d_{4,\frac 12}$}&1&3&
\centerline{$\mathcal J(e_0)=\alpha^3$,}
\centerline{ $\mathcal J(e_1)=\alpha^2$}&
\centerline{$\mathcal J(e_0)=\alpha^3$,}
\centerline{ $\mathcal J(e_1)=e_2$}&
\centerline{$\mathcal J(e_0)=\frac 12 e_2$,}
\centerline{ $\mathcal J(e_1)=e_3$}\\
\hline
\centerline{$\mathfrak d_{4,2}$}&1&3&
\centerline{$\mathcal J(e_0)=\alpha^3$,}
\centerline{$\mathcal J(e_1)=\alpha^2$}&
\centerline{$\mathcal J(e_0)=e_1$,}
\centerline{$\mathcal J(e_2)=\alpha^3$}&
\centerline{$\mathcal J(e_0)=e_1$,}
\centerline{$\mathcal J(e_2)=-\frac 12 e_3$}\\
\hline
\centerline{$\mathfrak d'_{4,0}$}&2&4&
\centerline{---}&
\centerline{$\mathcal J(e_0)=\alpha^3$,}
\centerline{$\mathcal J(e_1)=e_2$}&
\centerline{$\mathcal J(e_0)=e_3$,}
\centerline{$\mathcal J(e_1)=e_2$}\\
\hline
\centerline{$\mathfrak d'_{4,\lambda}$}
\centerline{$\lambda>0$}&1&1&
\centerline{$\mathcal J(e_0)=2\lambda\alpha^3$,}
\centerline{$\mathcal J(e_1)=\alpha^2$}&
\centerline{$\mathcal J(e_0)=\alpha^3$,}
\centerline{$\mathcal J(e_1)=e_2$}&
\centerline{$\mathcal J(e_0)=e_3$,}
\centerline{$\mathcal J(e_1)=e_2$}\\
\hline \caption{}\label{tabla}
\end{longtable}
\bigskip


\section{An example in dimension $6$}\label{six-dim}

In this section we exhibit an example of a six dimensional (non-nilpotent) solvable
Lie algebra $\mathfrak g_6$ admitting neither
symplectic nor complex  structures but having generalized complex
structures of types $1$ and $2$. This proves that Theorem $4.1$
fails for solvable Lie algebras of dimension six.
Examples of six dimensional nilpotent Lie algebras
having neither left invariant symplectic  nor
complex structures but with generalized complex structures
are given in~\cite{CG}.

Let us consider the solvable $6$-dimensional Lie algebra $\mathfrak
g_6$ defined by the structure equations
\begin{equation}\label{streq6}
d\alpha^i=0,\ (1\leq i\leq 4), \quad
d\alpha^5=\alpha^{12}+\mu\alpha^{15},\quad d\alpha^6=\alpha^{15}+\alpha^{34},\qquad
\mu\not=0.
\end{equation}
Let $\omega$ be a 2--form $\omega=\displaystyle\sum_{1\leq i<j\leq6} \omega_{ij}
\alpha^{ij}$. Then, one can check that
$\omega$ is closed if and only if
$$\omega_{16}=\omega_{26}=\omega_{36}=\omega_{46}=\omega_{56}=
\omega_{25}=\omega_{35}=\omega_{45}=0,$$
that is, $\omega$ is expressed as
\begin{equation}\label{simplectica6}
\omega=\sum_{i=2}^5 \omega_{1i}\;\alpha^{1i}+\sum_{i=3}^4 \omega_{2i}\;\alpha^{2i}+
\omega_{34}\;\alpha^{34}.
\end{equation}
But such a form $\omega$ is degenerate. This means that the Lie
algebra $\mathfrak g_6$ does not admit generalized complex
structures of type~0.

On the other hand, we consider the basis $\{X_i\}_{i=1}^6$ dual to
the basis of  1-forms $\{\alpha^i\}_{i=1}^6$. From the equations
\eqref{streq6} we get that the only non-zero Lie brackets are
\[
[X_1,X_2]=-X_5,\quad [X_1,X_5]=-X_6-\mu X_5,\quad [X_3,X_4]=-X_6.
\]
Let $J$ be an almost complex structure  on $\mathfrak g_6$ defined
by $J(X_i)=\displaystyle\sum_{j=1}^6a_{ij}X_j$.
Then, the components $(N_J)_{ij}^k$ of the Nijenhuis tensor
$N_J$ of $J$ satisfy
\[
(N_J)_{46}^3 =a_{63}{}^2,\quad
(N_J)_{36}^4 =-a_{64}{}^2,\quad
(N_J)_{56}^1 =\mu a_{51}a_{61}+a_{61}{}^2,\quad
(N_J)_{26}^1 =a_{51}a_{61}.
\]
Therefore, $\ a_{61}=a_{63}=a_{64}=0$ since $\ N_J=0$. Moreover,  $J^2=-\id$ implies
that
$a_{21}{}^2+a_{31}{}^2+a_{41}{}^2+a_{51}{}^2\not=0$, and the equations
\[
(N_J)_{26}^6 =a_{65}a_{21},\quad (N_J)_{36}^6 =a_{65}a_{31},\quad
(N_J)_{46}^6 =a_{65}a_{41},\quad (N_J)_{56}^6 =a_{65}a_{51},
\]
imply that $a_{65}=0$. Now, from
\[
(N_J)_{26}^5 =a_{62}a_{21},\quad (N_J)_{36}^5 =a_{62}a_{31},\quad
(N_J)_{46}^5 =a_{62}a_{41},\quad (N_J)_{56}^5 =a_{62}a_{51},
\]
we obtain $a_{62}=0$. But $-1=(J^2)_6^6= a_{66}{}^2$, which is not
possible. This proves that $\mathfrak g_6$ does not admit  complex
structures or, equivalently, it does not admit generalized complex
structures of type 3.

To describe generalized complex structures of types 1 or 2,  we take the
basis $\left\{X_i\right\}_{i=1}^{12}$ of $T^*\mathfrak g$ given by
\[
X_i=(X_i,0)\qquad \mbox{and} \qquad X_{i+6}=(0,\alpha^i),\qquad 1\leq i\leq 6.
\]
Notice that  the matrices $\mathcal J_i$ in \eqref{jotas4} have order $6$
\[
\mathcal J_1=\begin{pmatrix} a_{11}&a_{12}&a_{13}&a_{14}&a_{15}&a_{16}\\
 a_{21}&a_{22}&a_{23}&a_{24}&a_{25}&a_{26}\\
a_{31}&a_{32}&a_{33}&a_{34}&a_{35}&a_{36}\\
 a_{41}&a_{42}&a_{43}&a_{44}&a_{45}&a_{46}\\
a_{51}&a_{52}&a_{53}&a_{54}&a_{55}&a_{56}\\
 a_{61}&a_{62}&a_{63}&a_{64}&a_{65}&a_{66}\end{pmatrix},\qquad
\mathcal J_2=\begin{pmatrix} 0&b_{12}&b_{13}&b_{14}&b_{15}&b_{16}\\
-b_{12}&0&b_{23}&b_{24}&b_{25}&b_{26}\\-b_{13}&-b_{23}&0&b_{34}&b_{35}&b_{36}\\
-b_{14}&-b_{24}&-b_{34}&0&b_{45}&b_{46}\\-b_{15}&-b_{25}&-b_{35}&-b_{45}&0&b_{56}\\
-b_{16}&-b_{26}&-b_{36}&-b_{46}&-b_{56}&0\end{pmatrix}.
\]
Analogous calculations to those performed in the previous
sections allow us to obtain
\[
\begin{array}{lc}
a_{13}=a_{14}=a_{16}=a_{26}=a_{36}=a_{46}=0,&\\
b_{12}=b_{13}=b_{14}=b_{23}=b_{24}=b_{34}=0,&\\
b_{15}=\mu b_{16},\quad a_{15}=\mu a_{12}.
\end{array}
\]
Notice that $\det (\mathcal J_2)=0$, and so $\mathfrak g_6$ does not admit
generalized complex structures of type 0 as we already knew.

Let us consider
\[
\begin{array}{lc}
\left(\mathcal J^2\right)_1^1=-1=a_{11}{}^2+a_{12}(a_{21}+\mu a_{51})- b_{16}(c_{16}
+\mu c_{15}),&\\
\left(\mathcal J^2\right)_1^7=0=-2\mu^2 a_{12}b_{16}.
\end{array}
\]
This leads us to distinguish  the two following cases, that
exhaust all the possibilities:
\[
  \text{(i) } b_{16}\not=0\ \ \text{ but }\ \ a_{12}=0, \qquad \quad
  \text{ (ii) } b_{16}=0\ \ \text{ but }\ \ a_{12}\not=0.
\]
In case (i), assuming that $\rank \mathcal J_2=2$ and the integrabiblity of the
generalized almost complex structure, we get
$\left(\mathcal J^2\right)_5^5=(a_{55}-\mu a_{65})^2$, which is not
possible. So, a necessary condition for  having a generalized
complex structure is that $\rank \mathcal J_2=4$, that is,
\begin{equation}\label{elec1}
b_{25}\not=\mu b_{26}, \quad b_{35}\not=\mu b_{36}\quad\mbox{or}\quad b_{35}\not=\mu
b_{56}.
\end{equation}
Taking, for example, $b_{25}\not=\mu b_{26}$ we have the following
generalized complex structure of type $1$:
\[
\begin{array}{llll}
\mathcal J(X_1)=\alpha^6,&\mathcal J(X_2)=\alpha^5-\mu\alpha^6,&\mathcal
J(X_3)=X_4,&\mathcal J(X_4)=-X_3\\
\mathcal J(X_5)=-\alpha^2,&\mathcal J(X_6)=-\alpha^1+\mu\alpha^2,&\mathcal
J(\alpha^1)=\mu X_5+X_6,&\mathcal J(\alpha^2)=X_5,\\
\mathcal J(\alpha^3)=\alpha^4,&\mathcal J(\alpha^4)=-\alpha^3,&\mathcal
J(\alpha^5)=-\mu X_1-X_2,&\mathcal J(\alpha^6)=-X_1.
\end{array}
\]
Similar results are obtained for the remaining choices in
\eqref{elec1}.

In  case (ii), that is, $b_{16}=0$ and $a_{12}\not=0$, the condition
$\mathcal J^2=-\id$ implies that
\[
0=\left(\mathcal J^2\right)_1^8=-\mu a_{12}b_{25},\quad
0=\left(\mathcal J^2\right)_1^9=-\mu a_{12}b_{35},\quad
0=\left(\mathcal J^2\right)_1^{10}=-\mu a_{12}b_{45},
\]
\[
0=\left(\mathcal J^2\right)_1^{12}=a_{12}(b_{26}+\mu b_{56}),
\]
and then
\[
b_{25}=b_{35}=b_{45}=0, \qquad\mbox{and}\qquad b_{26}=-\mu b_{56}.
\]
Therefore, $\rank \mathcal J_2=0$ or $2$. If $\rank \mathcal J_2=0$,
we have $-1=\left(\mathcal J^2\right)_6^6= a_{66}{}^2$, which is not
possible. This means that $\mathfrak g_6$ does not admit generalized
complex structures of type 3, as we mentioned before. So  we must
have  $\rank \mathcal J_2=2$, that is,
\begin{equation}\label{elec2}
b_{36}\not=0,\qquad b_{46}\not=0\qquad\mbox{or}\qquad b_{56}\not=0.
\end{equation}
Consider $b_{56}\not=0$. Then for $b_{56}=1$ we obtain the
following generalized complex structure:
\begin{equation}\label{ecg6tipo2X}
\begin{array}{lll}
\mathcal J(X_1) =X_2,&\mathcal J(X_3) =X_4,&\mathcal J(X_5)=-\mu X_1-\alpha^6,\\
\mathcal J(X_6) =\alpha^5,&
\mathcal J(\alpha^1) =\alpha^2+\mu\alpha^5,&\mathcal J(\alpha^2)
=-\alpha^1+\mu X_6,\\ \mathcal J(\alpha^3) =\alpha^4,&\mathcal J(\alpha^6)
=-\mu X_2+X_5.
\end{array}
\end{equation}
Changing the basis $\{X_i,\alpha^i\}$ to  $\{Y_i,\beta^i\}$
defined by
\[
Y_i=X_i,\ \ i=1,2,3,4,6,\qquad \mbox{and} \qquad Y_5= X_5-\mu X_2,
\]
(and the corresponding change between the dual basis
$\{\alpha^i\}$ and $\{\beta^i\}$, resp.) the equations
\eqref{ecg6tipo2X} can be written, as stated in
Theorem~\ref{simple}, in the simplest form:
\[
\begin{array}{c}
\mathcal J(Y_1) =Y_2,\qquad \mathcal J(Y_3) =Y_4,\qquad \mathcal J(X_5) =-\beta^6,\\
\mathcal J(X_6) =\beta^5,\qquad \mathcal J(\beta^1) =\beta^2,\qquad \mathcal
J(\beta^3) =\beta^4,
\end{array}
\]
that is, the generalized complex structure $\mathcal J$ has type
2. For the remaining choices in \eqref{elec2} we obtain similar
results.

\smallskip

\noindent {\bf Acknowledgments.} We would like to thank the
referee for useful comments, which improved this paper.This work
has been partially
supported through grants MCyT (Spain) Project BFM2001-3778-C03-02,
MTM2005-08757-C04-02 and UPV 00127.310-E-15909/2004, ANPCyT
(Argentina) Project PICT3-13557 and Fundaci\'on Antorchas
(Argentina) Project 53900/3.


\vskip 1 cm

\noindent{\sf L.C. de Andr\'es and M. Fern\'andez:}
Departamento de Matem\'aticas, Facultad de Ciencia y Tecnolog\'{\i}a, Universidad del Pa\'{\i}s
Vasco, Apartado 644, 48080 Bilbao, Spain.

\noindent{\sl E-mail} de Andr\'es: {\tt luisc.deandres@ehu.es}\\
{\sl E-mail} Fern\'andez: {\tt marisa.fernandez@ehu.es}\\

\noindent{\sf M.L. Barberis and I. Dotti:} Facultad de Matem\'atica,
Astronom\'{\i}a y F\'{\i}sica,  Universidad Nacional de C\'ordoba, Ciudad Universitaria,
5000 C\'ordoba, Argentina.

\noindent {\sl E-mail} Barberis: {\tt barberis@mate.uncor.edu}\\
{\sl E-mail} Dotti: {\tt idotti@mate.uncor.edu}\\


\begin{thebibliography}{99} \frenchspacing

\bibitem{Ale}
D. Alekseevsky, J. Grabowsky, G. Marmo, P.W. Michor:
{\it Poisson structures on double Lie groups}, J. Geom. Phys. {\bf
26} (1998), 340--379.

\bibitem {ABDO}
 A. Andrada, M.L. Barberis, I. Dotti, G. Ovando: {\it Product
structures on four dimensional solvable Lie algebras}, Homology,
Homotopy  Appl.  {\bf 7} (1) (2005), 9--37,
arXiv:math.RA/0402234.

\bibitem {ABO}  A. Andrada, M.L. Barberis,  G. Ovando: {\it
Lie bialgebras of complex type and associated Poisson Lie groups},
preprint arXiv:math.DG/0610415.


\bibitem {BD}
 M.L. Barberis, I. Dotti: {\it Complex structures on affine motion
 groups}, Quart. J. Math. Oxford {\bf  55} (2004), 375--389.

\bibitem{BB}
O. Ben-Bassat: {\it Mirror symmetry and generalized complex
manifolds: Part I. The transform on vector bundles, spinors and
branes}, J. Geom. Phys. {\bf 56}  (4) (2006) 533--558,
arXiv:math.AG/0405303 v1.

\bibitem{BBB}
O. Ben-Bassat, M. Boyarchenko: {\it Submanifolds of generalized
complex manifolds}, J. Symplectic Geom. {\bf 2} (3) (2004),
309--355,
 arXiv:math.DG/0309013.

\bibitem{Cav}
G.R. Cavalcanti: {\it New  aspects of the $dd^c$-lemma}, Ph.D.
Thesis, University of Oxford, $2004$.

\bibitem {CG}
G.R. Cavalcanti, M. Gualtieri: {\it Generalized complex structures
on nilmanifolds},  J. Symplectic Geom. {\bf 2} (3) (2004), 393--410.

\bibitem{D} J. Dozias: {\it Sur les alg\`ebres de Lie r\'esolubles r\'eelles
de dimension inf\'erieure ou \'egale \`a $5$}, Th\`ese de 3 cycle,
novembre 1963, Facult\'e des Sciences de Paris.

\bibitem {FGG}
M.  Fern\'{a}ndez, M. Gotay, A. Gray:
{\it  Compact parallelizable four dimensional symplectic and complex
 manifolds}, Proc. Amer. Math. Soc. {\bf 103} (1988), 1209--1212.

\bibitem {Gualt}
M.  Gualtieri: {\it  Generalized complex geometry}, Ph.D. Thesis,
University of Oxford, $2003$; arXiv:math.DG/0401221.

\bibitem {Hitchin}
N. Hitchin: {\it  Generalized Calabi-Yau Manifolds}, Quart. J. Math. Oxford {\bf 54}
(2003), 281--308.

\bibitem {KN}
S. Kobayashi, K. Nomizu:  Foundations
of Differential Geometry vol. II, Interscience, 1969.

\bibitem {Mag}
L. Magnin:
{\it  Sur les algebres de Lie nilpotentes de dimension $\leq 7$}, J. Geom. Phys.
 {\bf 3} (1986), 119--131.

\bibitem{MR}
 A. Medina, P. Revoy: {\it  Groupes de Lie \`a structure symplectique invariante}, in
``Symplectic geometry, grupoids and integrable systems, S\'eminaire Sud-Rhodanien de
G\'eom\'etrie'' (P. Dazord et A. Weinstein Eds.), Mathematical Sciences Research
Institute publications, 247--266, New York-Berlin-Heidelberg, Springer Verlag 1991.

 \bibitem{Mu} G.M. Mubarakzyanov: {\it On solvable Lie algebras},
Izv. Vyssh. Uchebn. Zaved. Mat {\bf 32 } (1) (1963), 114--123
(Russian). Zbl 166,42.


\bibitem{O1} G. Ovando: {\it Invariant complex structures on
solvable real Lie groups}, Manuscripta Math. {\bf 103} (2000),
19--30.

\bibitem{O}
G. Ovando: {\it Complex, symplectic and K\"{a}hler structures on four
dimensional Lie groups}, Rev. Un. Mat. Argentina, {\bf 45} (2)
(2004), 55--67.

\bibitem{Z}  J. Patera, R.T. Sharp, P. Winternitz, H. Zassenhaus:
{\it Invariants of real low dimension Lie algebras}, J. Math.
Physics {\bf 17} (1976), 986--994.

\bibitem {Sal}
S. Salamon: {\it Complex structures on nilpotent Lie algebras},
 J. Pure Applied Algebra  {\bf 157} (2001), 311--333.

\bibitem {Sam}
H. Samelson: {\it A class of complex analytic manifolds},
 Portugal. Math. {\bf 12} (1953), 129--132.

\bibitem {Sa}
T. Sasaki: {\it Classification of left invariant
complex structures on $GL(2, \mathbb R)$ and $U(2)$}, Kumammoto J. Sci. (Math.)
 {\bf 14} (1981), 115--123.

 \bibitem{SJ} J.E. Snow: {\it Invariant Complex Structures
on Four Dimensional Solvable Real Lie Groups}, Manuscripta Math.
{\bf 66} \, (1990),  397--412.
\end{thebibliography}
\end{document}